\documentclass[12pt]{article}
\usepackage{amsfonts}

\def\Bbb{\mathbb}
\makeatletter

   \ifnum\@ptsize=0  \fi
   \ifnum\@ptsize=2  \fi
   \if@twoside\fi
   \parindent=\leftmargini

   \def\maketitle{\begin{center}\let\thanks=\footnote
      {\large\bf\@title}\par\bigskip\bigskip
      {\sc\@author}\par\bigskip
      {\rm\@date}
      \end{center}\bigskip\thispagestyle{empty}}

   \def\l@section{\@dottedtocline{1}{0em}{1.6em}}
   \def\l@subsection{\@dottedtocline{2}{1.5em}{3.2em}}

   \newlength{\additionalsectionskip} \additionalsectionskip=0cm
   \let\osection=\section
   \def\starsection*#1{\vspace{\additionalsectionskip}
      \osection*{\centering\normalsize\bf#1}\@afterindenttrue}
   \def\normalsection#1{\refstepcounter{section}\starsection*{\thesection. #1}
      \addcontentsline{toc}{section}{\protect\numberline{\thesection.}#1}}
   \def\section{\@ifnextchar*{\starsection}{\normalsection}}

   \def\thmindent{}
   \def\@begintheorem#1#2{\trivlist\item[\hskip\labelsep\thmindent{\bf#1 #2.}]\it}
   \def\@opargbegintheorem#1#2#3{%
      \trivlist\item[\hskip\labelsep\thmindent{\bf#1 #2 {\rm(#3)}.}]\it}

  \let\o@item=\@item
  \def\@item[#1]{\o@item[\rm #1]}

   \def\@biblabel#1{#1.}
   \let\othebibliography=\thebibliography
   \def\thebibliography#1{\small
      \def\@listi{\topsep=0cm\parsep=0cm\itemsep=0cm}\othebibliography{#1}}

   \newtheorem{satz}{Satz}[section]
   \makeatletter
   \@addtoreset{equation}{satz}

   \newtheorem{corollary}[satz]{Corollary}
   
   \newtheorem{example}[satz]{Example}
   \newtheorem{lemma}[satz]{Lemma}
   \newtheorem{proposition}[satz]{Proposition}
   \newtheorem{remark}[satz]{Remark}
   \newtheorem{remarks}[satz]{Remarks}
   \newtheorem{theorem}[satz]{Theorem}

   \newtheorem{thevarthm}[satz]{\varthmname}
   \newenvironment{varthm}[1]{\def\varthmname{#1}\begin{thevarthm}}
      {\end{thevarthm}\def\varthmname{}}
   \def\thmindent{}
   \newenvironment{varthm*}[1]{\trivlist\item[]{\thmindent\bf #1.}\it}
      {\endtrivlist}

   \def\epsilon{\varepsilon}                
   \def\phi{\varphi}                        
   \def\to{\longrightarrow}                 
   \def\mapsto{\mapstochar\longrightarrow}  
                   
   \def\bar{\overline}                      
   \def\hat{\widehat}                       
   \def\tilde{\widetilde}                   
   \def\({\left(}                           
   \def\){\right)}                          
   \def\O{{\cal O}}

   \def\set#1{\left\{\,#1\,\right\}}
   \def\be{\begin{eqnarray*}}
   \def\ee{\end{eqnarray*}}
   \def\matr#1#2#3#4{\left(\begin{array}{cc}#1&#2\\ #3&#4\end{array}\right)}
   \def\liste#1#2#3{\mbox{$#1_{#2},\dots,#1_{#3}$}} 
   \def\eqnref#1{(\ref{#1})}                        
   \def\Cal#1{{\cal#1}}

   \def\isom{\cong}                                 
   \def\isomto{\mathop{\longrightarrow}\limits^\sim}
   \def\with{\ \vrule\ }                                  
   
   \def\Bigwith{\with} 
   \def\Biggwith{\with}
   \def\tensor{\otimes}                             
   \def\rounddown#1{\left\lfloor#1\right\rfloor}    
   \def\inverse{^{-1}}                              

   \def\eqdef{=_{\rm def}}                          
   \def\operatorname#1{\mathop{\rm #1}\nolimits}
   \def\Bl{\operatorname{Bl}}
   \def\codim{\operatorname{codim}}
   \def\End{\operatorname{End}}
   \def\max{\operatorname{max}}
   \def\mult{\operatorname{mult}}
   \def\Nef{\operatorname{Nef}}
   \def\NS{\operatorname{NS}}
   \def\Pic{\operatorname{Pic}}
   \def\rank{\operatorname{rank}}
   \def\Sym{\operatorname{Sym}}

   \newenvironment{casearray}{%
      \left\{\begin{array}{cl}
   }{\end{array}\right.}

   \newenvironment{items}{\list{\labelitemi}{
      \parsep=0cm \itemsep=0cm \topsep=0cm \partopsep=\medskipamount
      \def\makelabel##1{\hss\llap{##1}}}
   }{\endlist}

   \newenvironment{pitems}{
      \def\item[##1]{\par\mbox{}\phantom{\rm(a)}\llap{\rm##1}}
   }{\par\noindent}

   \newenvironment{eqnarray+}{
      \begin{equation}\begin{array}{rcl}
   }{\end{array}\end{equation}}

   \def\lrskip{\parindent}
   \def\lreqn#1#2{
      \begin{trivlist}\item[]
         $\def\arraystretch{1.3}
         \hskip\lrskip\begin{array}{r}#1\end{array}$\hspace*{\fill} \\
         \hspace*{\fill}
         $\def\arraystretch{1.3}\begin{array}[t]{l}#2\end{array}\hskip\lrskip$
      \end{trivlist}
   }

   \def\proof{\trivlist\item[\hskip\labelsep\thmindent{\em Proof.}]}
   \def\proofof#1{\trivlist\item[\hskip\labelsep\thmindent{\em Proof of #1.}]}
   \def\endproof{\hspace*{\fill}\endproofsymbol\endtrivlist}
   \def\endproofsymbol{\frame{\rule[0pt]{0pt}{8pt}\rule[0pt]{8pt}{0pt}}}

   \renewenvironment{varthm*}[1]{\trivlist\item[]{\thmindent\bf #1.}\it}
      {\endtrivlist}

   \def\eps{\epsilon}
   \def\P{{\Bbb P}}
   \def\Q{{\Bbb Q}}
   \def\R{{\Bbb R}}
   \def\Z{{\Bbb Z}}
   \def\PP{{\Bbb P}^2}
   \def\PPP{{\Bbb P}^3}
   \def\and{\quad\mbox{ and }\quad}

   \def\N{N}
   \def\NE{N\hskip-0.1em E}
   \def\bNE{{}\hskip0.25em\bar{\hskip-0.25em N\hskip-0.1em E}}
   \def\NS{N\hskip-0.1em S}
   \def\NSR{\NS_{\Bbb R}}

   \def\Lsd{L_{\sqrt d}}
   \def\bilin#1#2{\left\langle#1,#2\right\rangle}

\begin{document}

   \abovedisplayshortskip=2pt plus 3pt
   \belowdisplayskip=14pt plus 3pt minus 7pt
   \belowdisplayshortskip=9pt plus 3pt minus 2pt
   \partopsep=7pt plus 2pt minus 0pt

   \title{Seshadri constants on algebraic surfaces}
   \author{Thomas Bauer%
   }
   \date{\phantom{July 27, 1998}}
   \maketitle
   \tableofcontents


\setcounter{section}{-1}
\section{Introduction}

   Seshadri constants are local invariants that are naturally associated
   to polarized varieties.  Except in the simplest cases they are
   very hard to control or to compute explicitly.
   The purpose of the present paper is to study these invariants on
   algebraic surfaces; we prove a number of explicit bounds on Seshadri
   constants and Seshadri sub-maximal curves, and we give complete results
   for abelian surfaces of Picard number one.

   In recent years
   there has been considerable interest in understanding
   the local positivity of ample line bundles on algebraic varieties.
   Motivated in part by the study of linear series in connection with
   Fujita's conjectures, Demailly \cite{Dem92}
   captured the concept of local positivity in the
   {\em Seshadri constant}, a real number $\eps(L,x)$ associated
   with an ample line bundle $L$ at a point $x$
   of an algebraic variety $X$, which in effect measures
   how much of the positivity of $L$ can be concentrated at $x$.
   Interest in Seshadri constant derives on the one hand from the fact
   that, via vanishing theorems,
   a lower bound on the Seshadri constants $\eps(L,x)$
   yields bounds
   on the number of points and jets that the adjoint series $\O_X(K_X+L)$
   separates.
   On the other hand it has become increasingly clear
   that Seshadri constants are highly interesting
   invariants of algebraic varieties quite in their own right.
   For instance, the papers
   \cite{Pao94} and \cite{Bau97a} address the question as to
   what kind of geometric information is encoded in them.
   We refer to Section \ref{section Seshadri} and to
   \cite[Section 1]{EinKueLaz95}
   for more on background and motivation.

   Even though on surfaces
   linear series are reasonably well understood
   thanks to powerful methods such as Reider's theorem, Seshadri constants
   are---as Demailly pointed out in \cite{Dem92}---extremely delicate
   already in the two-dimensional case.  For instance,
   if $X$ is a generic smooth surface of degree $d$ in $\P^3$, then
   the Seshadri constants of its hyperplane bundle are unknown when $d\ge 5$
   (cf.\ \cite{Bau97a} for the case $d=4$).
   In light of these facts it seems interesting to study Seshadri constants
   in the surface case, and to aim for explicit bounds or even
   explicit values.

   Our first result concerns
   surfaces that come with
   a fixed embedding into projective space.
   It is clear that in this case one has $\eps(L,x)\ge 1$ at all points
   $x$, and it is natural to ask under which circumstances
   equality holds and which small values bigger than $1$ can occur.
   Theorem \ref{theorem very ample} answers these questions.

   The second result deals with
   line bundles that are merely assumed
   to be ample.
   Work of Ein and Lazarsfeld \cite{EinLaz93b} on surfaces,
   and Ein-K\"uchle-Lazarsfeld \cite{EinKueLaz95}
   for the higher dimensional case,
   shows that
   there exist universal
   lower bounds on Seshadri constants
   if one restricts one's attention to very general points.
   Refinements of this type of results are due to
   K\"uchle and Steffens \cite{KueSte}, Steffens \cite{Ste}
   and Xu \cite{Xu95b}.
   On the other hand,
   well-known examples due to Miranda \cite[Proposition 5.12]{Laz93}
   show that there cannot exist
   universal lower bounds on Seshadri constants that are valid at
   {\em arbitrary} points.
   So 
   any bound on Seshadri constants at not necessarily very general
   points needs to take
   into account the geometry of the polarized surface $(X,L)$ in some way.
   We provide in Theorem \ref{theorem canonical slope} a bound in terms
   of the canonical slope of $L$, an invariant defined in terms of the
   nef cone of the surface.  We also carry out
   a closer analysis of Miranda's
   examples, which illustrates the interplay between this invariant and the
   Seshadri constant.

   We next study Seshadri sub-maximal curves at very general
   points, i.e.\ curves
   causing the Seshadri constant $\eps(L,x)$ to be below its maximal possible
   value $\sqrt{L^2}$ at these points.
   Naturally, it is of interest to find constraints
   on the existence of such curves, for instance upper bounds on their
   $L$-degree.  Now, as the results of Section \ref{section asf} will
   show, the curves computing Seshadri constants can have
   arbitrarily high degree, and---still worse---there cannot exist a
   bound involving only $L^2$.
   However,
   working variationally as in \cite{EinLaz93b}
   and \cite{Ste},
   we show in Theorem \ref{theorem very general}
   that there does exist an explicit bound for curves
   causing $\eps(L,x)$ to be less than $\sqrt{L^2}-\delta$,
   for given $\delta>0$,
   in terms of $L^2$ and $\delta$.
   In a similar vein, we ask in Section \ref{section number sub-maximal}
   for the number of sub-maximal curves at a fixed (not necessarily
   very general) point, and we provide an explicit upper bound on this
   number in \ref{proposition number}.

   In Sections \ref{section asf} to \ref{section multiple} we
   focus on the study of Seshadri constants on abelian surfaces.
   Our motivation to investigate this class of surfaces is twofold:
   First, there are up to now hardly any non-trivial examples where
   Seshadri constants are explicitly known.
   And secondly, the study of Seshadri constants on abelian varieties
   in general
   has just recently gained considerable attention
   (see \cite{Nak96}, \cite{Laz97} and \cite{Bau2}).
   Notably, in \cite{Laz97} Lazarsfeld has established a
   surprising connection between the Seshadri constant $\eps(X,L)$
   of an abelian variety $(X,L)$ and
   a metric invariant, the
   minimal period length $m(X,L)$.
   This allows on the one hand to get lower bounds on $\eps(X,L)$
   as soon as a lower bound on $m(X,L)$ is available, as is the case
   for principal polarizations by work of Buser and Sarnak
   \cite{BusSar94} and for polarizations of arbitrary type
   by \cite{Bau2}. On the other hand, it can be used to show that
   certain abelian varieties (such as Jacobians \cite{Laz97} and Prym
   varieties \cite{Bau2}) have a period of unusually short length.
   For the surface case, previous joint work \cite{BauSze}
   of the author with Szemberg
   gave an upper bound on the Seshadri constant $\eps(X,L)$ of an
   abelian surface $(X,L)$ involving the solutions of a diophantine
   equation. In some
   cases, this upper bound could be shown to be equal to $\eps(X,L)$,
   leading to the speculation that this might always be true
   if $(X,L)$ is generic.
   In Theorem \ref{theorem asf} we will complete the picture by showing
   that this is in fact the case.
   A nice feature of this result is that it allows to explicitly
   compute the Seshadri constants---as well as
   the unique irreducible curve that accounts for it---for
   a whole class of surfaces.
   It also
   shows that Seshadri constants have an intriguing number-theoretic
   flavor in this case.
   In view of the papers \cite{Laz97} and \cite{Bau2}
   it would be most interesting to know if a similar formula exists
   for abelian varieties in any dimension.

   In Section 7 we classify
   the nef cones of abelian surfaces, which we think is interesting
   both in view of applications to
   Seshadri constants and as a complement to \cite{Bau1}.
   Finally we consider Seshadri constants along
   finite sets
   in Section 8. This generalized notion,
   which appears implicitly already
   in Nagata's famous conjecture \cite{Nag59},
   has been studied previously
   by Xu (see \cite{Xu94} and \cite[5.16]{Laz93}) and K\"uchle
   \cite{Kue96}.  Our purpose here is to provide
   a lower bound for multiple-point
   Seshadri constants on abelian surfaces and to study their relationship
   with the one-point constant.

\begin{varthm*}{\it Notation and Conventions}
\rm
   We work throughout over the field $\Bbb C$ of complex numbers.
   The symbol $\equiv$ denotes numerical equivalence of divisors or
   line bundles, whereas linear equivalence will be denoted by $\sim$.
   For a real number $x$ we denote by $\rounddown x$ is
   its round-down (integer part).
   We will say that a property holds for a {\em very general} point
   of a variety $X$, if it holds off the union of countably many proper
   closed subvarieties of $X$.
\end{varthm*}


\section{Seshadri constants}\label{section Seshadri}

   We briefly recall in this section the definition of Seshadri constants
   and their relationship with linear series.  For a more detailed
   exposition we refer to
   \cite[Section 5]{Laz93} and \cite[Section 1]{EinKueLaz95}.

   Consider a smooth projective variety $X$ and a nef line bundle
   $L$ on $X$.
   Let $$f:Y=\Bl_x(X)\to X$$ be the blow-up of $X$ at a point $x\in X$
   and $E=f\inverse(x)$ the exceptional divisor.
   The {\em Seshadri constant} $\eps(L,x)$ is by definition the
   maximal real number $\eps$ such that
   the line bundle $f^*L-\eps E$ is nef, i.e.\
   $$
      \eps(L,x)\eqdef\sup\set{\eps\in\R\with f^*L-\eps E\mbox{ nef}} \ .
   $$
   It is elementary that this can be equivalently expressed on the variety
   $X$ itself
   as the infimum
   $$
      \eps(L,x)=\inf\set{\frac{L\cdot C}{\mult_x C}\with
      C\subset X\mbox{ irreducible curve through } x } \ .
   $$
   Intuitively, the number $\eps(L,x)$ measures how much of the positivity
   of $L$ can be concentrated at the point $x$.
   The name Seshadri constant derives from the fact that
   by Seshadri's criterion, $L$ is ample if and only if its
   {\em global Seshadri constant}
   $$
      \eps(L)=\inf\set{\eps(L,x)\with x\in X}
   $$
   is positive.

   There are interesting characterizations of Seshadri constants in terms
   of linear series. Recall that, for an integer $s\ge 0$,
   a line bundle $B$ on $X$ (or the linear series $|B|$)
   is said to {\em separate
   $s$-jets at $x$}, if $B$ admits global sections with arbitrarily
   prescribed $s$-jet at $x$, i.e.\ if
   the evaluation map
   $$
      H^0(X,B)\to H^0(X,B\tensor\O_X/\Cal I_x^{s+1})
   $$
   is onto.
   Let $s(B,x)$ denote the largest integer $s$ such that $|B|$ separates
   $s$-jets at $x$.
   The relationship of the separation of jets with
   Seshadri constants can then be summarized as follows:

\begin{proposition}\label{Seshadri and jets}
   Let $X$ be a smooth projective variety, $x\in X$ a point,
   and $L$ an ample line bundle on $X$.
   \begin{pitems}
   \item[(a)]
      For every line bundle $M$ on $X$, one has
      $$
         \eps(L,x)=\limsup_{k\to\infty}\frac{s(M+kL, x)}k \ .
      $$
   \item[(b)]
      If $k$ and $s$ are non-negative integers satisfying the inequality
      $$
         k>\frac{s+\dim(X)}{\eps(L,x)} \ ,
      $$
      then
      the adjoint linear series $|K_X+kL|$ separates $s$-jets at $x$.
   \end{pitems}
\end{proposition}

   When $M=\O_X$ then statement (a) is Theorem 6.4 in \cite{Dem92};
   the fact that is holds for arbitrary $M$ can be shown using
   arguments as in the proof of \cite[Proposition 1.1(b)]{EinKueLaz95}.
   Part (b) is an application of
   Kawamata-Viehweg vanishing on the blow-up of $X$
   (see \cite[Proposition 1.1(a)]{EinKueLaz95}).
   Note a subtle but crucial difference between (a) and (b):
   Whereas statement (a) tells us that $\eps(L,x)$ is determined by
   the asymptotic behaviour of the series $|M+kL|$ for $k\gg 0$, where
   $M$ is {\em any} given line bundle (for instance $M=\O_X$ or $M=K_X$),
   the statement in (b), which gives information
   about 
   the series $|K_X+kL|$ for $k\ge k_0$ with an explicit value for $k_0$,
   does not remain true
   in general when $K_X$ is replaced by an arbitrary
   line bundle $M$ (e.g.\ by $\O_X$).


\section{Very ample line bundles}
\label{section very ample}

   Consider a smooth projective surface $X$ and a very ample line bundle
   $L$ on $X$. If $C$ is an irreducible curve on $X$ and $x\in C$ a
   point, then there is certainly
   a divisor $D\in|L|$ passing through $x$ and
   meeting $C$ properly. Thus the Seshadri constant $\eps(L,x)$, and
   hence the global Seshadri constant $\eps(L)$, is always at least one.
   This elementary argument already gives
   the best possible lower bound in general, since
   if the surface $X$ contains a line (when embedded via the linear
   series $|L|$), then equality $\eps(L)=1$ is attained.
   It is then natural to ask whether this is the only case where
   $\eps(L)=1$ occurs, and what the next possible values of $\eps(L)$
   for a very ample line bundle are. Theorem \ref{theorem very ample}
   below answers these questions. When dealing with a very ample line
   bundle $L$ on $X$, we will
   identify $X$ with the surface in $\P(H^0(X,L))$ obtained by
   embedding $X$ via the linear series $|L|$ and we will
   write
   $$
      \eps(X, x) \eqdef \eps(\O_X(1), x) = \eps(L, x) \ .
   $$
   and
   $$
      \eps(X) \eqdef \eps(\O_X(1)) = \eps(L) \ .
   $$
   We refer to these numbers simply as {\em the Seshadri constants of
   $X$}, tacitly suppressing the (fixed) choice of the projective embedding.

\begin{theorem}\label{theorem very ample}
\begin{pitems}
\item[(a)]
   Let $X\subset\P^N$ be a smooth irreducible surface. Then
   $\eps(X)=1$ if and only if $X$ contains a line.

\item[(b)]
   For $d\ge 4$ let $\Cal S_{d,N}$
   denote the space of smooth irreducible
   surfaces of degree $d$ in $\P^N$ that do not contain any lines.
   Then
   $$
      \min\set{\eps(X)\with X\in\Cal S_{d,N}} = \frac d{d-1} \ .
   $$

\item[(c)]
   If $X$ is a surface in $\Cal S_{d,N}$ and $x\in X$ is a point such that
   the local Seshadri constant $\eps(X,x)$ satisfies the inequalities
   $1<\eps(X,x)<2$,
   then it is of the form
   $$
      \eps(X,x) = \frac ab \ ,
   $$
   where $a,b$ are integers with $3\le a\le d$ and $a/2<b<a$.

\item[(d)]
   All rational numbers $a/b$ with $3\le a\le d$ and $a/2<b<a$
   occur as local Seshadri constants of smooth irreducible surfaces
   in $\P^3$ of degree $d$.
\end{pitems}
\end{theorem}

\proof
   (a)
   The ''if'' part being obvious we assume $\eps(X)=1$ and show that
   $X$ contains a line.  Our first claim is then:
   \begin{equation}\label{line claim}
      \parbox{0.8\textwidth}{
         There is an irreducible curve $C\subset X$ and a point $x\in X$
         such that $\deg(C)=\mult_x(C)$.
       } 
   \end{equation}
   To prove \eqnref{line claim} we assume by way of contradiction that
   there is a sequence $(C_n,x_n)_{n\ge 0}$ of irreducible curves
   $C_n\subset X$ and points $x_n\in X$ such that
   $$
      \eps_n \eqdef \frac{\deg(C_n)}{\mult_{x_n}(C_n)} \to 1
      \quad\mbox{, but $\eps_n>1$ for all $n\ge 0$.}
   $$
   We may choose the $C_n$ such that $\eps_n<2$ for all $n\ge 0$.
   Since $h^0(X,\O_X(1))\ge 4$, there is for every $n\ge 0$ a divisor
   $D_n\in|\O_X(1)|$ with $\mult_{x_n}(D_n)\ge 2$.  If $D_n$ and $C_n$
   were to meet properly, then we would have
   \be
      \deg(C_n)&=&D_nC_n \\
         &\ge& \mult_{x_n}(D_n)\cdot\mult_{x_n}(C_n) \\
         &\ge& 2\mult_{x_n}(C_n) \ ,
   \ee
   and hence $\eps_n\ge 2$, which contradicts our assumption on $\eps_n$ above.
   So $C_n$ must be a component of $D_n$, thus we get the estimate
   \begin{equation}\label{degree estimate}
      \deg(C_n)\le D_n^2 = \deg(X) \ .
   \end{equation}
   We may assume that the curves $C_n$ are chosen in such a way
   that $\eps_n\le 1+\frac 1n$, hence
   $$
      \mult_{x_n}(C_n) < \deg(C_n) \le \(1+\frac 1n\)\mult_{x_n}(C_n) \ .
   $$
   The fact that $\deg(C_n)$ is an integer then implies that
   $\mult_{x_n}(C_n)\ge n$.  So we get $\deg(C_n)>n$, which
   contradicts the upper bound \eqnref{degree estimate} on the degree
   of the $C_n$.  This establishes the claim \eqnref{line claim}.

   To prove the statement in the theorem,
   we are now going to show that the curve
   $C$ in \eqnref{line claim} is in fact a line.  To this end,
   consider the projection
   $$
      \pi:\P^N-\{x\} \to \P^{N-1}
   $$
   from $x$
   onto a hyperplane $\P^{N-1}\subset\P^N$.
   Let $H\subset\P^N$ be a hyperplane passing through $x$ and not
   containing $C$. Then
   the equality $\deg(C)=\mult_x(C)$
   implies that the intersection $H\cap C$ is supported entirely on
   the point $x$.  So we find that
   the image $\pi(C-\{x\})$ does not meet the
   generic hyperplane $H$.  But then $\pi(C-\{x\})$ must be finite, and
   hence $C$ is a line.

   (b)
   We first show that $\eps(X)\ge\frac d{d-1}$ for all $X\in\Cal S_{d,N}$.
   So let $X\in\Cal S_{d,N}$. We may assume $\eps(X)<2$. Then there is an
   irreducible curve $C\subset X$ and a point $x\in X$ such that
   \begin{equation}\label{less than 2}
      \frac{L\cdot C}{\mult_x C} < 2
   \end{equation}
   Let $H\subset\P^N$ be a hyperplane containing the tangent plane
   $T_xX$. Then \eqnref{less than 2} implies that $C$ is a component
   of the hyperplane section $X\cap H\in|L|$. Since this holds for all
   hyperplanes $H$ containing $T_xX$, we conclude that $C\subset X\cap
   T_xX$, i.e.\ $C$ lies in the tangent plane at $x$. But then we know
   that
   $$
      \mult_xC\le L\cdot C-1 \ ,
   $$
   since $C$ cannot be a line, and therefore
   $$
      \frac{L\cdot C}{\mult_xC}
      \ge \frac{L\cdot C}{L\cdot C-1}\ge \frac d{d-1} \ .
   $$
   This shows that $\eps(X)\ge d/(d-1)$, as claimed.
   The fact that the minimum is taken on by some surface is a special
   case of statement (d) which we will prove below.

   (c)
   Suppose $C\subset X$ is an irreducible curve and $x\in X$ a point
   such that
   $$
      1<\frac{L\cdot C}{\mult_xC}<2 \ .
   $$
   It follows as in (b) that $C$ is a
   component of a divisor $D$, whose support is contained in the
   tangent plane $T_xX$, so that we have
   $$
      1<\deg C\le d \and \frac{\deg C}2<\mult_xC<\deg C \ .
   $$
   The intersection $X\cap T_xX$ consists
   of only finitely many irreducible curves, and the Seshadri
   constant $\eps(X,x)$ is computed by one of these curves. This shows
   that $\eps(X,x)$ is of the form $a/b$ as claimed.

   (d)
   It remains to show that all rational numbers $a/b$, where
   $a$ and $b$ satisfy the conditions
   specified in the theorem, occur in this way. To this end, given $a$
   and $b$, we choose according to Lemma \ref{lemma curves}(a) below an
   irreducible curve $C_0\subset\PP$ of degree $a$ with a point $x$ of
   multiplicity $b$. Further, we take a smooth curve $C_1\subset\PP$
   of degree $d-a$
   not passing through $x$. By statement (b) of Lemma \ref{lemma
   curves}, there is a smooth surface $X\subset\PPP$ such that
   the divisor $C_0+C_1$ is a hyperplane section of $X$. According to
   the arguments in the proof of (b), the curve $C$ computing
   $\eps(X,x)$ is a component of the intersection $X\cap T_xX$,
   and therefore $C=C_0$. So we conclude
   $$
      \eps(X,x)=\frac{L\cdot C_0}{\mult_xC_0}=\frac ab \ .
   $$
   This completes the proof of the theorem.
\endproof

\begin{lemma}\label{lemma curves}
\begin{pitems}
\item[(a)]
   Let $p\in\PP$ be a point and let $m$ and $d$ be
   positive integers with $m<d$. Then there are
   irreducible curves of degree $d$ with a point of
   multiplicity $m$ at $p$.

\item[(b)]
   For every reduced divisor $D\subset\PP$ there is a smooth
   surface $X\subset\PPP$ such that $D$ is a hyperplane section of
   $X$.
\end{pitems}
\end{lemma}

   Note that in (b) the assumption that $D$ be reduced is essential: A
   non-reduced divisor can never be a hyperplane section of a smooth
   surface in projective three-space, since the
   Gau{\ss} map of a smooth hypersurface of degree $\ge 2$ is finite.

\proof
   (a)
   Fix a point $p\in\PP$ and consider on the blow-up $f:X\to\PP$ at $p$
   with exceptional divisor $E$ over $p$ the line bundle
   $$
      M_{d,m}\eqdef f^*\O_{\PP}(d)-mE \ .
   $$
   Since $\O_{\PP}(d)$ is $d$-jet ample, the line bundle
   $M_{d,m}$ is globally generated whenever $m\le d$ (cf.\
   \cite[Lemma 3.1]{BelSom93b}). Moreover,
   for $m<d$ we have $M_{d,m}^2>0$, so that the linear series
   $|M_{d,m}|$ is not composed with a pencil. Now take an irreducible
   element $C_0\in|M_{d,m}|$. Its direct image $C=f_*C_0$ satisfies
   our requirements.

   (b)
   Let $d$ denote the degree of $D$. We may certainly assume $d\ge 2$.
   Since $D$ is reduced, we can choose a curve $D'\subset\PP$ of
   degree $d-1$ meeting $D$ transversely. Let $g$ and $g'$ be affine
   equations of $D$ and $D'$ respectively, and consider the surface
   $X\subset\PPP$ with the affine equation
   $$
      f(x_1,x_2,x_3)=g(x_1,x_2)+x_3\cdot g'(x_1,x_2) \ .
   $$
   It contains the divisor $D$ as its intersection with the plane
   $H=\{x_3=0\}$, so it remains to show that $X$ is smooth. Suppose
   then to the contrary that there is a singularity $p=(p_1,p_2,p_3)$
   of $X$. Looking at the derivative of $f$ with respect to the
   coordinate $x_3$ we obtain $g'(p_1,p_2)=0$, so that the equation
   $f(p)=0$ implies $g(p_1,p_2)=0$. This means that the point
   $p'=(p_1,p_2)$ is a point of intersection of the divisors
   $D$ and $D'$.
   Therefore the 1-jet $j^1_p(f)$ of $f$ at $p$ is given by
   $$
      0=j^1_p(f)=j^1_{p'}(g)+j^1_p(x_3\cdot g')
         =j^1_{p'}(g) + p_3\cdot j^1_{p'}(g') \ .
   $$
   But this says that the divisors $D$ and $D'$ have a common tangent at
   $p'$, contradicting the choice of $D'$.
\endproof

\begin{remark}\rm
   Consider a smooth surface $X\subset\P^3$ of degree $\ge 3$.
   By Theorem \ref{theorem very ample}(a), the locus
   \begin{equation}\label{locus of lines}
      \set{x\in X\with\eps(X,x)=1}
   \end{equation}
   is the union of all lines on $X$. While this locus
   is empty for generic $X$, special surfaces often contain quite a
   number of lines. In any event, there are always only finitely many
   lines on a {\em smooth} surface, and
   it is an interesting, yet unsolved, classical problem
   to determine the {\em maximal} number
   $\ell(d)$ of lines that can lie on a smooth
   surface of degree $d$ in $\P^3$ for any given $d\ge 3$.
   The only numbers that are explicitly known are $\ell(4)=64$
   (see \cite{Seg43}), and of course $\ell(3)=27$.
   For $d\ge 5$, the best general bounds available at present
   are
   $$
      3d^2 \le \ell(d) \le 11d^2-28d-8
   $$
   (see \cite[Sect.\ 5.1]{CapHarMaz95} for the lower bound
   and \cite[Sect.\ 4]{Seg43}
   for the upper bound, cf.\ also \cite[Sect.\ 2.4 and 5]{Miy84}).
   It would be interesting to know if the characterization
   \eqnref{locus of lines} could be used to derive upper bounds on $\ell(d)$.
   More generally, one may ask if there are explicit
   bounds on the degree of the loci
   $$
      \set{x\in X\with\eps(X,x)\le a} \qquad\mbox{ for } a\ge 1 \ .
   $$
\end{remark}


\section{Bounds on global Seshadri constants}
\label{sect canonical slope}

   In Section \ref{section very ample} we dealt with the Seshadri
   constants of very ample line bundles. For these line bundles the
   lower bound $\eps(L)\ge 1$ is obvious. If we consider ample
   bundles, however, the situation changes quite dramatically: One
   knows from Miranda's examples (see \cite[Proposition 5.12]{Laz93})
   that the global Seshadri constant
   $\eps(L)$ can become arbitrarily small. More precisely, there
   are sequences of polarized varieties $(X_k,L_k)$ such that
   $\eps(L_k)< 1/k$.
   So any lower bound on $\eps(L)$ has to involve the
   geometry of the polarized variety $(X,L)$ in some way. A closer
   analysis of Miranda's examples
   (see the proof of Proposition \ref{proposition Miranda} below)
   shows that the positivity of the line
   bundles $L_k$ with respect to the canonical divisor $K_X$ get
   smaller and smaller as $k$ grows. Theorem \ref{theorem canonical
   slope} will show that this behaviour is necessary for the fact that the
   Seshadri constant gets small.
   The number that accounts for the relation of $L$ and $K_X$ in this
   context is the
   slope of $L$ relative to $K_X$ with respect to the nef cone $\Nef(X)$
   in the vector space $N^1(X)$ of
   real-valued classes of codimension one on $X$, i.e.\
   the minimum
   $$
      \sigma(L)\eqdef\min\set{s\in\R\with\O_X(sL-K_X)\mbox{ is nef}} \
      .
   $$
   We will simply refer to the number $\sigma(L)$ as the {\em canonical
   slope} of $L$. The intuition here is that a
   line bundle $L$ with large canonical slope is 'bad' in the sense that
   one needs a high multiple of $L$ to reach the positivity of $K_X$.

   Using this notion, our result can be stated as follows:

\begin{theorem}\label{theorem canonical slope}
   Let $X$ be a smooth projective surface and $L$ an ample line bundle
   on $X$. Then the global Seshadri constant of $L$ is bounded in terms
   of the canonical slope of $L$ by
   $$
      \eps(L)\ge\frac 2{1+\sqrt{4\sigma(L)+13}} \ .
   $$
\end{theorem}

\begin{remarks}\rm
   (a) For $(X,L)=(\P^2,\O_{\P^2}(1))$ one has $\sigma(L)=-3$ and
   $\eps(X)=1$, so that equality holds in Theorem
   \ref{theorem canonical slope}.

   (b) As for another example, consider a K3 surface $X$. For any
   ample line bundle $L$ on $X$ one has $\sigma(L)=0$, hence the
   theorem gives
   $$
      \eps(L)\ge\frac 2{1+\sqrt{13}}=0.434\dots
   $$
   It is easy to see that the
   optimal statement in this situation is $\eps(L)\ge 1/2$.

   (c) Suppose that $X$ is a surface of general type with $K_X$ ample.
   We have $\sigma(K_X)=1$, and hence
   $$
      \eps(K_X)\ge\frac 2{1+\sqrt{17}}=0.390...
   $$
   Of course, we do not expect that this particular bound is sharp.
   It would however
   be interesting to know how far it is from being optimal.

   (d) On any smooth projective surface one has the lower bound
   $\sigma(L)\ge -3$ (for
   instance by Mori's theorem or by Reider's criterion for global
   generation of line bundles).
   On the other hand, the invariant
   $\sigma(L)$ can become arbitrarily large, as already
   smooth surfaces in $\P^3$ show (see (e)).

   (e)
   One certainly cannot expect
   that $\sigma(L)$ alone fully accounts for the
   behaviour of the Seshadri constant. Consider for instance a smooth
   surface $X\subset\P^3$ of degree $d$ and take $L=\O_X(1)$. Then one
   always has $\sigma(L)=d-4$, whereas the value of $\eps(L)$
   depends on the geometry of $X$ (see e.g.\cite{Bau97a}).

\end{remarks}

\proofof{Theorem \ref{theorem canonical slope}}
   Let $C\subset X$ be an irreducible curve, $x\in X$ a point, and
   $m=\mult_x C$. The idea is to use that fact that a point of
   multiplicity $m$ causes the
   geometric genus of a curve
   to drop by at least ${m \choose 2}$, in order to
   derive an upper bound on $m$.
   Specifically, we have by the adjunction formula
   $$
      1 + \frac12C(C+K_X)=p_a(C)\ge p_a(C)-p_g(C)\ge{m\choose 2} \ ,
   $$
   so that we get
   $$
      m(m-1)\le 2+C(C+K_X) \ .
   $$
   Now, by assumption, $\O_X(\sigma(L)L-K_X)$ is nef, so that in particular
   $K_X\cdot C\le\sigma(L)L\cdot C$,
   and therefore
   $$
      C(C+K_X)\le C^2+\sigma(L)L\cdot C \ .
   $$
   This gives an upper bound on the multiplicity $m$ of $C$ at $x$:
   $$
      m(m-1)\le C^2+\sigma(L)L\cdot C+2 \ ,
   $$
   which, upon using the Hodge index theorem, implies
   $$
      m\le\frac12+\sqrt{\frac{(L\cdot C)^2}{L^2}+\sigma(L)L\cdot C+\frac94} \ .
   $$
   The upshot of these considerations is that
   $$
      \eps(L)\ge\min_{d\ge 1}\frac d{\frac12+
      \sqrt{\frac{d^2}{L^2}+\sigma(L)d+\frac94}} \ .
   $$
   One checks now, for instance using a little elementary analysis,
   that this minimum is taken on at $d=1$.
\endproof

   In order to shed some more light on the interplay of the invariants
   $\eps(L)$ and $\sigma(L)$, we will now have a closer look at Miranda's
   examples. In fact, we will give a generalized version of
   \cite[Proposition 5.12]{Laz93}, which shows
   that the occurence of small
   Seshadri constants is not as exceptional as it might appear judging
   from Miranda's examples, which are originally constructed as
   certain blow-ups of $\P^2$. The
   following proposition shows that actually a suitable blow-up of {\em any}
   surface with Picard number one contains ample line bundles with
   Seshadri constants below any prescribed bound.

\begin{proposition}\label{proposition Miranda}
   Let $X$ be a smooth projective surface of Picard number one. Then:
   \begin{pitems}
   \item[(a)]
      For every integer $r>0$ there are line bundles $L_k$, $k\ge 1$, on
      suitable blow-ups $Y_k$ of $X$ such that
      $L_k^2\ge r^2$,
      $L_k\cdot C\ge r$ for all curves $C$ on $Y$, but
      $$
         \eps(L_k)\le\frac 1k \ .
      $$
   \item[(b)]
      For $X=\P^2$ and $r=1$ the bundles $L_k$ can be chosen in such a
      way that for every real number $\delta>0$ we have
      $$
         \eps(L_k)-\eps_0(L_k)<\frac{\delta}{k} \qquad\mbox{ for } k\gg 0 \ ,
      $$
      where $\eps_0(L_k)$ denotes the lower bound (in terms of the canonical
      slope $\sigma(L_k)$) given by Theorem \ref{theorem canonical slope}.
   \end{pitems}
\end{proposition}

\begin{remark}\rm
   Note in particular that, while $\eps(L_k)$ gets arbitrarily small,
   the stated intersection
   properties of $L_k$ allow one to
   achieve that 
   the adjoint linear series $|K_X+L_k|$ generates an arbitrarily
   prescribed number
   of jets at any point, and hence has an arbitrarily large global
   Seshadri constant $\eps(K_X+L_k)$.
   More concretely, given any integer $s>0$, it suffices to choose
   $r\ge s^2+4s+5$ in order to force the linear series $|K_X+L_k|$
   to generate $s$-jets at any point of $X$
   (e.g.\ by \cite[Corollary 7.5]{Laz93}).
   But this implies by Proposition \ref{Seshadri and jets}
   that $\eps(K_X+L_k)\ge s$.
   So, from the point of view of Seshadri constants,
   the linear series $|L_k|$ itself does not
   benefit from the numerical positivity of $L_k$ (expressed in terms of
   its self-intersection and its intersection with curves), whereas the
   adjoint linear series $|K_X+L_k|$ does.
\end{remark}

   In order to prove the proposition, we will need the following
   lemma:

\begin{lemma}\label{lemma red}
   Let $L$ be an ample line bundle whose class generates the N\'eron-Severi
   group of a smooth surface $X$. Consider for $d>0$ the space
   $\Cal R_d$ of reducible divisors in the linear series $|\O_X(dL)|$. Then
   there is a constant $c$ such that
   $$
      \codim(\Cal R_d, |\O_X(dL)|)\ge dL^2+c \ .
   $$
\end{lemma}

\proofof{Lemma \ref{lemma red}}
   By Serre vanishing, there is an integer $n_0$ such that for $n\ge
   n_0$ one has
   \begin{equation}\label{Serre for nL}
      H^i(X,\O_X(nL))=0 \quad\mbox{ for }i>0 \ .
   \end{equation}
   Moreover, the integer $n_0$ can be chosen in such a way that the
   vanishing remains true even if $L$ is replaced by a numerically
   equivalent line bundle (cf.\ \cite[Theorem 5.1]{Fuj83}).
   To prove the lemma,
   it is enough to show that for $1\le a<d$ and for every
   $P\in\Pic^0(X)$
   \begin{equation}\label{to prove}
      \codim\(|\O_X(aL)\tensor P|+|\O_X((d-a)L)\tensor P^{-1}|, |\O_X(dL)|\)
      \ge dL^2+c
   \end{equation}
   for some constant $c$ independent of $d$ and $P$,
   and it suffices to consider $d\ge 2n_0$.
   We distinguish between two cases. Suppose first that $a\ge n_0$ and
   $d-a\ge n_0$. Thanks to \eqnref{Serre for nL} we have then by Riemann-Roch
   \lreqn{
      h^0(X,\O_X(dL))
      -h^0(X,\O_X(aL)\tensor P)-h^0(X,\O_X((d-a)L)\tensor P\inverse)
   }{
      =\chi(X,\O_X(dL))
         -\chi(X,\O_X(aL)\tensor P)-\chi(X,\O_X((d-a)L)\tensor P\inverse) \\
      =a(d-a)L^2-\chi(O_X) \\
      \ge (d-1)L^2-\chi(\O_X) \\
      = dL^2+\mbox{const}
      \ ,
   }
   which proves the assertion \eqnref{to prove} in this case.
   In the alternative case, we may by symmetry assume that
   $a<n_0$ and $d-a>n_0$. Then we have
   \begin{equation}\label{d-a ineq}
      h^0(X,\O_X((d-a)L)\tensor P\inverse)\le h^0(X,\O_X((d-1)L)\tensor P\inverse)
      \ ,
   \end{equation}
   and, using the abbreviation
   $$
      b \eqdef\max\set{h^0(X,\O_X(kL)\tensor P)\with 1\le k\le n_0,\
      P\in\Pic^0(X)} \ ,
   $$
   we obtain upon using \eqnref{Serre for nL} and \eqnref{d-a ineq}
   the estimate
   \lreqn{
      h^0(X,\O_X(dL))
      -h^0(X,\O_X(aL)\tensor P)-h^0(X,\O_X((d-a)L)\tensor P\inverse)
   }{
      \ge h^0(X,\O_X(dL))
         -b-h^0(X,\O_X((d-1)L)\tensor P\inverse) \\
      = \chi(X,\O_X(dL))
         -b-\chi(X,\O_X((d-1)L)\tensor P\inverse) \\
      = \frac12(2d-1)L^2-b-\frac12 L\cdot K_X \\
      = dL^2+\mbox{const} \ ,
   }
   and this completes the proof of the lemma.
\endproof

\proofof{Proposition \ref{proposition Miranda}}
   (a) We will follow Miranda's construction to exhibit line bundles
   with the properties asserted in the proposition.
   We start by choosing an
   ample generator $H$ of $\NS(X)$, fixing an integer $k\ge 1$
   and setting $m=rk$.
   For sufficiently large $d>0$, the line bundle $\O_X(dH)$ will be
   $(m+1)$-jet ample, so that we can find an irreducible curve
   $D\in|dH|$ with a point $x$ of multiplicity $\ge m$.
   In view of Lemma \ref{lemma red} we can arrange, by possibly
   increasing $d$, that there is a pencil $P$ of irreducible curves in
   $|dH|$ containing $D$.  We may further assume that $P$ has $d^2H^2$
   distinct base points $\liste p1{d^2H^2}\in X$. Consider the blow-up
   $f:Y=Y_k\to X$ of $X$ at these points and the induced pencil
   $$
      \hat P=f^*P-\sum_{i=1}^{d^2H^2} E_i \ ,
   $$
   where $E_i=f^{-1}(p_i)$. The proper transform $\hat D$ of $D$ has
   multiplicity $\ge m$ at the point $f^{-1}(x)$.
   Fix now an integer $a\ge 2$ and consider the line bundle
   $$
      L=L_k\eqdef\O_Y(r(a\hat D+E_1)) \ .
   $$
   For the intersection numbers of $L$ we have the bounds
   \begin{equation}\label{L self}
      L^2=r^2(2a-1)\ge r^2
   \end{equation}
   and
   \begin{eqnarray+}\label{L curves}
      L\cdot\hat D&=&r(a\hat D^2+E_1\cdot \hat D)=r>0 \\
      L\cdot E_1&=&r(a\hat D\cdot E_1+E_1^2)=r(a-1)>0 \ .
   \end{eqnarray+}%
   The map $Y\to\P^1$ induced by $\hat P$ is a fibration with
   irreducible fibres and with section $E_1$. Therefore, by the
   Nakai-Moishezon criterion, the inequalities
   \eqnref{L self} and \eqnref{L curves}
   imply that $L$ is ample.
   Due to the existence of the singular curve $\hat D$,
   its Seshadri constant is bounded from
   above by
   $$
      \eps(L)\le\eps(L,f^{-1}(x))\le\frac{L\cdot\hat D}{\mult_{f^{-1}(x)}\hat D}
      \le\frac rm=\frac 1k \ ,
   $$
   whereas of course
   $L\cdot C\ge r$ for all curves $C$ on $Y$,
   since $L$ is an $r$-th power.

   (b) Note first that for $X=\P^2$, $r=1$ and $d\gg 0$
   we may take $d=m+1=k+1$. We now
   determine an upper bound on
   the canonical slope of the bundles $L_k$ for $k\gg 0$.
   Writing $E=\sum_{i=1}^{d^2}E_i$, we have
   $$
      K_{Y}=f^*K_X+E=-3f^*H+E \ ,
   $$
   so that we find
   $$
      K_{Y}\cdot\hat D=d(d-3),\
      K_{Y}\cdot E_1=-1,\
      K_{Y}\cdot L_k=ad^2-3ad-1,\
   $$
   and $K_{Y}^2=9-d^2$.  The line bundle $s L_k-K_{Y}$
   therefore satisfies
   \be
      (s L_k- K_{Y})\hat D&=&s-d(d-3) \\
      (s L_k- K_{Y})E_1&=&s(a-1)+1
   \ee
   and
   $$
       (s L_k- K_{Y})^2=(2a-1)s^2-2s(ad^2-3ad-1)+9-d^2 \ .
   $$
   Fix now a real number $\eta>1$.
   One checks then that for $d\gg 0$ the Nakai-Moishezon criterion
   implies that $s L_k-K_Y$ is ample for $s\ge\eta d^2$, and hence
   $$
      \sigma(L_k)<\eta(k+1)^2 \quad\mbox{ for } k\gg 0 \ .
   $$
   We therefore get the estimate
   $$
      \eps(L_k)-\eps_0(L_k)\le\frac 1k-\frac 2{1+\sqrt{4\eta(k+1)^2+13}}
      \ ,
   $$
   and for $k\gg 0$ the latter expression gets smaller than $\delta/k$,
   if $\eta$ is chosen sufficiently close to $1$.
   This completes the proof of the proposition.
\endproof


\section{The degree of sub-maximal curves}

   In this section we show how the techniques from \cite{EinLaz93b}
   can be used to derive an explicit bound on the degrees of
   the irreducible curves leading to sub-maximal Seshadri constants
   at very general points.
   Specifically,
   suppose that $C$ is an irreducible curve and $x\in X$ a point such that
   the quotient
   $$
      \eps_{C,x} \eqdef \frac{L\cdot C}{\mult_x C}
   $$
   is less than $\sqrt{L^2}$. In other words, the curve $C$ causes the
   Seshadri constant $\eps(L,x)$ to be at most $\eps_{C,x}<\sqrt{L^2}$.
   We will briefly refer to curves with this property as {\em Seshadri
   sub-maximal curves}.
   From the point of view of Seshadri constants, these are the most
   interesting curves on $X$, because they account for the failure of
   $L$ to have maximal positivity.
   It is therefore highly desirable to obtain as much information about them
   as possible.
   The following result provides an upper bound on the degree of
   a Seshadri sub-maximal curve $C$ in terms of $\eps_{C,x}$.

\begin{theorem}\label{theorem very general}
   Let $X$ be a smooth projective surface and let $L$ be an ample line
   bundle on $X$. Further, let $x\in X$ be a very general point and
   $C\subset X$ an irreducible curve passing through $x$ such that
   $$
      \eps_{C,x} < \sqrt{L^2} \ .
   $$
   Then the degree of $C$ with respect to $L$ is bounded as follows:
   $$
      L\cdot C < \frac{L^2}{\sqrt{L^2}-\eps_{C,x}} \ .
   $$
\end{theorem}

   So, roughly speaking, the theorem says that
   only curves of small
   degree can force $\eps(L, x)$ to be small at very general points.

\begin{remark}\rm
   It is also useful to think of the theorem as giving an upper bound on the
   self-intersection of $C$. In fact, combining the inequality in the
   theorem with the Hodge index theorem yields the bound
   $$
      C^2 < \frac{L^2}{(\sqrt{L^2}-\eps_{C,x})^2} \ .
   $$
   Consider for instance the case when $\eps(L,x)\le \sqrt{L^2}-1$.
   The theorem implies then the existence of an irreducible curve
   $C$ passing through $x$ such that
   $$
      L\cdot C<L^2 \quad\mbox{ and }\quad C^2<L^2 \ .
   $$
\end{remark}

   In the proof of the theorem we
   will make use of the following result of Ein and Lazarsfeld in the
   spirit of \cite{Ste}.

\begin{proposition}[Ein-Lazarsfeld \cite{EinLaz93b}]\label{proposition EL}
   Let $X$ be a smooth projective surface and let
   $(C_t)_{t\in T}$ be a non-trivial 1-parameter family of irreducible
   curves $C_t\subset X$.  Suppose that $(x_t)_{t\in T}$ is
   a family of points $x_t\in C_t$
   and $m$ an integer such that
   $$
      \mult_{x_t}C_t \ge m
   $$
   for all $t\in T$. Then
   $$
      C^2_t \ge m(m-1) \ .
   $$
\end{proposition}

\proofof{Theorem \ref{theorem very general}}
   Let $m=\mult_x C$. Since $x$ is very general in $X$, there exists a
   non-trivial family $(C_t)_{t\in T}$ of irreducible curves $C_t\subset X$
   and a family $(x_t)_{t\in T}$ of points $x_t\in C_t$
   such that $\mult_{x_t}C_t\ge m$ and $(C_{t_0},x_{t_0})=(C,x)$ for some
   $t_0\in T$.
   Proposition \ref{proposition EL} then implies in particular
   \begin{equation}\label{EL bound}
      C^2\ge m(m-1) \ .
   \end{equation}
   Suppose now that $\alpha$ is a real number with
   $$
      \frac{L\cdot C}{m}<\alpha\le\sqrt{L^2}
   $$
   From these inequalities we obtain
   $
      \alpha L\cdot C<\alpha^2m\le mL^2
   $
   and hence
   \begin{equation}\label{alpha bound}
      \alpha\cdot\frac{L\cdot C}{L^2} < m \ .
   \end{equation}

   Now assume by way of contradiction that
   $$
      \frac{L^2}{L\cdot C} \le \sqrt{L^2}-\frac{L\cdot C}m \ .
   $$
   This implies that some multiple of the rational number $L^2/L\cdot C$
   is contained in the interval $(L\cdot C/m, \sqrt{L^2}]$, say  
   $$
      \frac{L\cdot C}m < k\frac{L^2}{L\cdot C} \le \sqrt{L^2}
   $$
   with a suitable integer $k$.
   Taking $\alpha=kL^2/L\cdot C$ and using the inequality
   \eqnref{alpha bound}
   we then have
   $$
      k < m \ .
   $$
   The crucial point of the proof is now an elementary diophantine argument
   in the spirit of \cite{Ste}: Since $k$ is an integer, the previous
   inequality implies $k\le m-1$. This slight improvement on the bound
   suffices to establish a contradiction. In fact,
   combining the inequality $k\le m-1$ with
   the bound \eqnref{EL bound} and with the Hodge index theorem, one obtains
   \be
      m(m-1)&\le& C^2 \\
         &\le& \sqrt{\frac{C^2}{L^2}}\cdot L\cdot C \\
         &<& m\alpha\sqrt{\frac{C^2}{L^2}} \\
         &=& mk \frac{L^2}{L\cdot C}\sqrt{\frac{C^2}{L^2}} \\
         &\le& mk \\
         &\le& m(m-1) ,
   \ee
   which is absurd, and this completes the proof of the theorem.
\endproof

   As an application we give a quick proof in the surface case of a
   result by Nakamaye \cite{Nak96},
   which characterizes the abelian surfaces of
   Seshadri constant one.

\begin{corollary}
   Let $(X,L)$ be a polarized abelian surface with $\eps(L)=1$. Then
   $(X,L)$ is a polarized product of elliptic curves,
   $$
      X=E_1\times E_2, \quad L=\O_X(d(E_1\times 0)+(0\times E_2)) \ ,
   $$
   where $d=L^2/2$.
\end{corollary}

\proof
   Fix a point $x\in X$. By assumption, there is
   for every $\delta>0$ an irreducible
   curve $C\subset X$ such that
   $$
      \frac{L\cdot C}{\mult_xC}<1+\delta \ .
   $$
   Since on abelian varieties the Seshadri constant is independent of
   the point, we can apply the
   theorem to get the inequality
   \begin{equation}\label{bound elliptic}
      C^2<\frac{L^2}{\(\sqrt{L^2}-(1+\delta)\)^2}
   \end{equation}
   If $L^2\ge 6$, then for small $\delta$
   the value of the expression on the right hand
   side is less than 2, so that we find $C^2=0$. Thus $C$ is an elliptic
   curve, and hence $\mult_xC=1$. This implies $L\cdot C=1$ and the
   assertion follows immediately.
   If $L^2\le 4$, then the inequality \eqnref{bound elliptic} implies
   $C^2\le 2$. So either $C^2=0$, where we conclude as before, or else
   $C^2=2$. In the latter case $C$ is a hyperelliptic curve
   of genus 2, and therefore again
   $\mult_xC=1$. This implies $L\cdot C=1$, which
   however is impossible by the Hodge index theorem.
\endproof


\section{On the number of sub-maximal curves}\label{sect lemma xi}
\label{section number sub-maximal}

   The canonical slope of an ample line bundle, which was used in Sect.\
   \ref{sect canonical slope} to obtain a lower bound on the global
   Seshadri constant, will also come into play when we consider
   the following
   enumerative question:  Let $L$ be an ample line bundle on a smooth
   projective surface $X$ and a point $x\in X$. Given a 
   real number $a>0$, what can we say about the number
   $\nu(L,x,a)$ of irreducible curves such that
   $$
      \frac{L\cdot C}{\mult_x C} < a \ ,
   $$
   if it is finite at all?
   Of course, one cannot expect to
   actually determine the number $\nu(L,x,a)$
   in general. However, we can give an explicit upper bound in terms of
   the canonical slope $\sigma(L)$ when $a$ is a rational number
   $<\sqrt{L^2}$:

\begin{proposition}\label{proposition number}
   Let $X$ be a smooth projective surface and $L$ an ample line bundle
   on $X$.  Suppose that a point $x\in X$ and a rational number
   $a<\sqrt{L^2}$ is given.  Set
   $$
      \delta = (\sigma(L)\cdot L^2+a)^2-8\cdot\chi(\O_X)(L^2-a^2) \ .
   $$
   Then
   $$
      \nu(L,x,a) \le kL^2 \ ,
   $$
   where $k>\sigma(L)$ is an integer such that the number $k\cdot a$
   is integral
   and, in case $\delta\ge 0$, such that
   $$
      k>\frac{\sigma(L)\cdot L^2+a+\sqrt{\delta}}{2(L^2-a)} \ .
   $$
\end{proposition}

   The idea for the proof of the proposition lies in the following
   useful observation:

\begin{lemma}\label{lemma xi}
   Let $X$ be a smooth projective surface and $L$ an ample line bundle
   on $X$.  Given a real number $\xi>0$ and a point $x\in X$,
   suppose that for some $k>0$ there is a divisor $D\in |\O_X(kL)|$
   such that
   $$
      \frac{L\cdot D}{\mult_x D}\le \xi\sqrt{L^2} \ .
   $$
   Then every irreducible curve $C\subset X$ satisfying the inequality
   $$
      \frac{L\cdot C}{\mult_x C}< \frac 1{\xi}\sqrt{L^2}
   $$
   is a component of $D$.
\end{lemma}

   For instance,
   this implies that an irreducible curve $C\in|\O_X(kL)|$
   with $L\cdot C/\mult_x C<\sqrt{L^2}$ computes the Seshadri constant
   $\eps(L,x)$.

\proof
   Suppose to the contrary that $D$ and $C$ intersect properly.  Then
   we get
   \be
      kL\cdot C&=&D\cdot C\ge \mult_x D\cdot\mult_x C \\
      &>& \frac{L\cdot D}{\xi\sqrt{L^2}} \cdot
      \frac{\xi L\cdot C}{\sqrt{L^2}} \\
      &=& kL\cdot C \ ,
   \ee
   and this is a contradiction.
\endproof

\proofof{Proposition \ref{proposition number}}
   The idea is simple: Find a divisor $D\in|\O_X(kL)|$ such that
   the quotient $L\cdot D/\mult_x D$ is sufficiently small;
   its degree will then
   by means of Lemma \ref{lemma xi} give an upper bound for
   $\nu(L,x,a)$.
   Turning to the details, let
   $\sigma=\sigma(L)$ and $k>\sigma$. We then have
   $$
      H^i(X,\O_X(kL))=H^i(X,\O_X(K_X+(k-\sigma)L+(\sigma L-K_X)))=0
      \quad\mbox{ for } i > 0
   $$
   by Kodaira vanishing, since $(k-\sigma)L$ is ample and $\sigma L-K_X$
   is nef.  Therefore
   \be
      h^0(X,\O_X(kL))&=&\chi(X,\O_X(kL)) \\
      &=&\chi(\O_X)+\frac 12 kL(kL-K_X) \\
      &=&\chi(\O_X)+\frac 12 kL\((k-\sigma)L+(\sigma L-K_X)\) \\
      &\ge&\chi(\O_X)+ \frac 12 k(k-\sigma)L^2 \ .
   \ee
   On the other hand, we have for $m>0$
   $$
      h^0(X,\O_X/\Cal I_x^m)={m+1\choose 2} \ ,
   $$
   so that the linear series $|\O_X(kL)\tensor\Cal I_x^m|$ will be
   non-empty as soon as
   $$
      \chi(\O_X)+\frac 12 k(k-\sigma)L^2-\frac 12 m(m+1) > 0 \ .
   $$
   If we take $m=k\cdot a$, which by assumption is an integer, then
   this condition
   is equivalent to the quadratic inequality
   \begin{equation}\label{quadr ineq}
      k^2(L^2-a^2)-k(\sigma L^2+a)+2\chi(\O_X) > 0 \ .
   \end{equation}
   So if $\delta$, its discriminant, is negative, then
   $|\O_X(kL)\tensor\Cal I_x^{ka}|\ne\emptyset$,
   since by assumption $a<\sqrt{L^2}$.
   If $\delta\ge 0$, then the linear series in question will be
   non-empty whenever $k>k_0$, where $k_0$ is the bigger root of
   the quadratic polynomial in
   \eqnref{quadr ineq}.

   In either case, taking a divisor $D\in|\O_X(kL)\tensor\Cal I_x^{ka}|$,
   we have
   $$
      \frac{L\cdot D}{\mult_x D}=\frac{kL^2}{\mult_x D}
      \le\frac{kL^2}{ka}=\xi\sqrt{L^2} \ ,
   $$
   where we set $\xi\eqdef\sqrt{L^2}/a$.
   If now $C\subset X$ is an irreducible curve with
   $$
      \frac{L\cdot C}{\mult_x C} < a = \frac 1{\xi}\sqrt{L^2} \ ,
   $$
   then, by Lemma \ref{lemma xi}, it is a component of $D$.
   This implies the assertion.
\endproof


\section{Seshadri constants of polarized abelian surfaces}
\label{section asf}

   Consider a polarized abelian variety $(X,L)$.
   By homogeneity, the Seshadri constant
   $\eps(L,x)$ is independent of the point $x\in X$, so it is an invariant
   of the polarized variety $(X,L)$.
   We will denote it by $\eps(X,L)$.
   There has been recent interest in the study of
   Seshadri constants of abelian varieties:
   Using symplectic blowing up in the spirit of
   \cite{McDPol94},
   Lazarsfeld has established an interesting connection between Seshadri
   constants and minimal period lengths, leading in particular
   to a lower bound on
   $\eps(X,L)$ for the principally polarized case.
   Generalizing the approach of Buser and Sarnak in \cite{BusSar94},
   a lower bound on $\eps(X,L)$ for arbitrary polarizations has been given
   in \cite{Bau2}.

   For the surface case, where one hopes for more
   specific results, an upper bound on $\eps(X,L)$ involving the solutions
   of a diophantine equation
   was given in \cite{BauSze}.  An interesting consequence of this result
   is that on abelian surfaces Seshadri constants are always rational.
   In certain cases the upper bound was shown to be equal to $\eps(X,L)$,
   and it was tempting to hope that
   this might always be true
   if $(X,L)$ is general.
   In Theorem \ref{theorem asf} we will complete the picture by showing
   that this is in fact the case.
   A nice feature of this result is that it allows to explicitly
   compute the Seshadri constants for a whole class of surfaces.
   It also allows to determine the unique irreducible curve that
   computes $\eps(X,L)$.

   We show:

\begin{theorem}\label{theorem asf}
   Let $(X,L)$ be a polarized abelian surface of type $(1,d)$, $d\ge 1$,
   such that $\NS(X)\isom\Z$.
   \begin{pitems}
   \item[(a)]
      If $\sqrt{2d}$ is rational, then $\eps(X,L)=\sqrt{2d}$.
   \item[(b)]
      If $\sqrt{2d}$ is irrational, then
      $$
         \eps(X,L)=2d\cdot\frac{k_0}{l_0}=\sqrt{2d\(1-\frac 1{\ell_0^2}\)}
         =\frac{2d}{\sqrt{2d+\frac 1{k_0^2}}}<\sqrt{2d} \ ,
      $$
      where $(k_0,\ell_0)$ is the primitive solution of Pell's equation
      $$
         \ell^2-2d k^2=1 \ .
      $$

      There is (up to translation) a unique irreducible curve $C\subset X$
      such that
      $$
         \eps(X,L)=\frac{L\cdot C}{\mult_x C}\quad\mbox{ for some } x\in C \ ,
      $$
      and we have either
      $$
         \O_X(C)\equiv\O_X(k_0L)\quad\mbox{ and }\quad\mult_x C=\ell_0
      $$
      or
      $$
         \O_X(C)\equiv\O_X(2k_0L)\quad\mbox{ and }\quad\mult_x C=2\ell_0 \ .
      $$
      Moreover, the point $x$ is the only singularity of the curve $C$.
   \end{pitems}
\end{theorem}

   In the proof
   we will apply the following useful lemma, which follows
   from \ref{lemma xi}:

\begin{lemma}\label{useful lemma}
   Let $X$ be a smooth projective surface, $x\in X$ a point, and $L$
   an ample line bundle on $X$. If there is a divisor $D\in|\O_X(kL)|$ for
   some $k>0$ satisfying
   $$
      \frac{L\cdot D}{\mult_x D}<\sqrt{L^2} \ ,
   $$
   then every irreducible curve $C\subset X$ with
   $$
      \frac{L\cdot C}{\mult_x C}<\sqrt{L^2}
   $$
   is a component of $D$.
\end{lemma}

\proofof{Theorem \ref{theorem asf}}
   Assertion (a) follows from Steffens' result \cite[Proposition 1]{Ste}
   to the effect that on any surface of Picard number one we have
   the lower bound
   $$
      \eps(L,x)\ge\rounddown{\sqrt{L^2}}
   $$
   for very general $x\in X$.

   As for (b):
   Replacing $L$ by a suitable algebraically equivalent line bundle,
   we may assume to begin with that $L$ is symmetric.
   It was shown in \cite{BauSze} that the linear series
   $|\O_X(2k_0L)|$ then contains an even symmetric divisor $D$ such that
   $\mult_{e_1} D\ge 2\ell_0$, where $e_1$ is a fixed halfperiod on $X$.
   Thus, as in \cite{BauSze}, we have the upper bound
   $$
      \eps(X,L)\le\frac{L\cdot D}{\mult_{e_1}D}\le\frac{k_0}{\ell_0}L^2 \ .
   $$
   Suppose now by way of contradiction that there is an irreducible
   curve $C\subset X$ passing through $e_1$ such that
   \begin{equation}\label{sesh for C}
      \frac{L\cdot C}{\mult_{e_1}C}<\frac{k_0}{\ell_0}L^2 \ .
   \end{equation}
   Let $\iota:X\to X$ denote the $(-1)$-involution on $X$.
   We have $\mult_{e_1}\iota^*C=\mult_{e_1}C$ and
   $L\cdot\iota^*C=L\cdot C$.
   Therefore,
   since $C$ is algebraically equivalent to a multiple of $L$, Lemma
   \ref{useful lemma} applies to $C$ and shows that the curves $C$ and
   $\iota^*C$ coincide, i.e.\ that $C$ is symmetric.
   Lemma \ref{useful lemma} also implies that $C$ appears as a component
   of $D$, hence we have
\begin{equation}\label{C data}
   \O_X(C)\equiv\O_X(k_1L) \quad\mbox{ with } k_1\le 2k_0
   \mbox{ and } m_1\eqdef\mult_{e_1}C\le\mult_{e_1}D \ .
\end{equation}

   Consider the blow-up $f:\tilde X\to X$ of $X$ at the sixteen halfperiods
   $\liste e1{16}$ and the projection $\pi:\tilde X\to K$ onto the smooth
   Kummer surface $K$ of $X$.
   Since $C$ is symmetric, its
   proper transform
   $$
      C'=f^*C-\sum_{i=1}^{16}\mult_{e_i}C\cdot E_i
   $$
   descends to an irreducible curve $\bar C\subset K$.
   We claim that
   \begin{equation}\label{1 on Kummer}
      h^0(K,\O_K(\bar C))=1 \ .
   \end{equation}
   In fact, if the linear series $|\O_K(\bar C)|$ were to contain a pencil,
   then this would give us a pencil of curves in $|\O_X(C)|$ with the same
   multiplicities at halfperiods as $C$. In particular, we would then
   have infinitely many irreducible curves satisfying \eqnref{sesh for C},
   which however by Lemma \ref{useful lemma} is impossible.
   This establishes \eqnref{1 on Kummer}.

   Now, \eqnref{1 on Kummer} implies $(\bar C)^2=-2$,
   since the exact sequence
   $$
      0\to\O_K(-\bar C)\to\O_K\to\O_{\bar C}\to 0
   $$
   tells us that
   $H^i(K,\O_K(\bar C))=0$ for $i>0$, so that
   by Riemann-Roch
   $$
      1=h^0(K,\O_K(\bar C))=\chi(K,\O_K(\bar C))=2+\frac 12(\bar C)^2 \ .
   $$
   We conclude that, using the abbreviation $m_i=\mult_{e_i}C$,
   \begin{equation}\label{mult proper}
      k_1^2\cdot 2d-\sum_{i=1}^{16} m_i^2
      =C^2-\sum_{i=1}^{16} m_i^2=(C')^2=(\pi^*\bar C)^2
      =2(\bar C)^2=-4 \ ,
   \end{equation}
   so that we obtain the lower bound
   \begin{equation}\label{lower for C}
      k_1^2\cdot 2d-m_1^2\ge -4 \ .
   \end{equation}
   On the other hand, we have the upper bound
   \begin{equation}\label{upper for C}
      k_1^2\cdot 2d-m_1^2<0 \ ,
   \end{equation}
   since in the alternative case the inequality
   $k_1/m_1\ge1/\sqrt{2d}$ would imply
   $$
      \frac{L\cdot C}{m_1}=\frac{k_1}{m_1}L^2\ge\sqrt{2d} \ ,
   $$
   a contradiction with \eqnref{sesh for C}.

   So, by \eqnref{lower for C} and \eqnref{upper for C}, there are
   only four possible values for the difference $k_1^2\cdot 2d-m_1^2$.
   We will deal with these cases separately.

   {\em Case 1.} Suppose that $k_1^2\cdot 2d-m_1^2=-4$.
   Then $m_1$ is necessarily an even number.
   From \eqnref{mult proper} we see that $m_i=0$ for $i>1$,
   so that in particular the multiplicities of $C$ at all halfperiods are
   even.  This implies that the symmetric line bundle $\O_X(C)$ is
   totally symmetric and is therefore algebraically equivalent
   to an even multiple of $L$.
   So $k_1$ is even as well.  But then the pair $(k_1/2, m_1/2)$
   is a solution of Pell's equation $\ell^2-2d k^2=1$.
   By the minimality of the solution $(k_0,\ell_0)$ we then have
   $k_1\ge 2k_0$ and $m_1\ge 2\ell_0$, and consequently by \eqnref{C data}
   $$
      k_1=2k_0 \and D=C \ .
   $$
   But this of course makes \eqnref{sesh for C} impossible.

   {\em Case 2.} Suppose that $k_1^2\cdot 2d-m_1^2=-3$.
   In this case $m_1$ is an odd number and, looking at the equation modulo
   4, we see that $k_1$ must be odd as well.  The symmetric line bundle
   $\O_X(C)$ then has $q$ odd halfperiods, where $q\in\{4,6,10,12\}$
   (cf.\ \cite[Section 5]{BirLan90}).  But we see from \eqnref{mult proper}
   that $m_1^2-\sum_{i=1}^{16} m_i^2=-1$, so that $C$ passes
   through only one halfperiod apart from $e_1$, which implies $q=2$,
   a contradiction.

   {\em Case 3.} Suppose that $k_1^2\cdot 2d-m_1^2=-2$.
   This is similar to the previous case:  Now $m_1$ is even and $k_1$
   is odd, as we see again by looking at the equation modulo 4.
   It follows from \eqnref{mult proper} that $C$ passes through only two
   halfperiods apart from $e_1$, and we get the same kind of contradiction
   as in Case 2.

   {\em Case 4.} Finally, suppose that  $k_1^2\cdot 2d-m_1^2=-1$.
   In this case the pair $(k_1,m_1)$ solves Pell's equation
   $\ell^2-2d k^2=1$, and the minimality of the solution $(k_0,\ell_0)$
   implies
   $$
      k_1=k_0 \and D=2C \ ,
   $$
   which does not allow \eqnref{sesh for C}.

   We now show the assertions about $C$.
   First, the uniqueness of $C$ is clear from Lemma \ref{useful lemma}.
   Further,
   since $\eps(X,L)$ is computed by $C$, and since
   $L\cdot D/\mult_{e_i}D=\eps(X,L)$, we must
   have $D=a\cdot C$ for some integer
   $a\ge 1$.
   The proof so far
   shows that either $D=2C$ (corresponding
   to Case 4) or $D=C$ (corresponding to Case 1).
   It remains to show that $e_1$ is the only singular point of $C$.
   The adjunction formula on $K$ tells us that
   $$
      p_a(\bar C)=1+\frac 12(\bar C)^2=0 \ ,
   $$
   so that in any event $C$ is smooth outside of the sixteen halfperiods.
   Further, we have either
   $$
      \sum_{i>1} m_i^2=0 \quad\mbox{ or }\quad \sum_{i>1} m_i^2=3 \ ,
   $$
   and this shows that $e_1$
   is the only halfperiod at which $C$ is singular.

   This completes the proof of the theorem.
\endproof

\begin{remark}\rm
   The statement about the numerical equivalence class of $C$ in part
   (b) of the theorem leaves two possibilities: either $C\equiv k_0L$
   or $C\equiv 2k_0L$.
   Let us stress here that
   both cases actually occur:
   If $2d+1$ is a square, then $(k_0,\ell_0)=(1,\sqrt{2d+1})$ is
   the minimal solution of Pell's equation, and
   the proof of \cite[Theorem A.1(c)]{BauSze} shows that $C\equiv k_0L$
   in this case.
   On the other hand, for $d=1$ we have
   $C\equiv 2k_0L$;
   this follows from the fact that
   $(k_0,\ell_0)=(2,3)$ and
   that
   the image of the unique divisor $\Theta\in|L|$ under
   the multiplication map $X\to X$, $x\mapsto 2x$, is an {\em irreducible}
   curve in $|4L|$ with multiplicity 6 at the origin (cf.\ \cite{Ste}).
\end{remark}

   The theorem implies in particular that $\eps(X,L)$ can be arbitrarily
   close to $\sqrt{L^2}$. Furthermore, it implies
   that the degree of the curve computing $\eps(X,L)$
   can be arbitrarily large, and that it cannot be bounded in terms of
   $L^2$ only:

\begin{corollary}\label{cor close}
   For every real number $\delta>0$ and every integer
   $N>0$ there is an integer
   $d>0$ such that for every polarized abelian surface $(X,L)$
   of type $(1,d)$ with $\NS(X)\isom\Z$ the following conditions hold:
   \begin{items}
   \item[(a)]
      $\sqrt{L^2}-\eps(X,L)<\delta$.
   \item[(b)]
      The unique irreducible curve $C\subset X$ that computes $\eps(X,L)$
      at $x\in X$ satisfies the inequalities
      $$
         L\cdot C>N\cdot L^2\quad\mbox{ and }\quad
         \mult_x C>N\sqrt{L^2} \ .
      $$
   \end{items}
\end{corollary}

\proof
   This follows essentially from the fact that for suitable $d$
   the solutions of Pell's equation are arbitrarily large.
   Specifically,
   let $d\ge 1$ be an integer such that $\sqrt{2d}$ is irrational, and let
   $p_n/q_n$, $n\ge 0$, be the convergents of $\sqrt{2d}$.
   One knows that
   for every solution $(k,\ell)$ of Pell's equation
   $\ell^2-2d k^2=1$,
   the rational number
   $\ell/k$ is one of the convergents $p_n/q_n$ (see e.g.\
   \cite[Chapter 10]{Hua81}).
   The sequences $(p_n)$ and $(q_n)$ have the following properties:
   $$
      p_{n+1}>p_n\ , q_{n+1}>q_n\ , p_0=a_0, q_0=1\ , p_1=a_1a_0+1\ ,
      q_1=a_1 \ ,
   $$
   where
   $$
      a_0=\rounddown{\sqrt{2d}}\quad\mbox{ and }\quad
      a_1=\rounddown{\frac1{\sqrt{2d}-a_0}} \ .
   $$
   We certainly have $(\ell,k)\ne(p_0,q_0)$, so that for the minimal
   solution $(k_0,\ell_0)$ we have the
   lower bound
   $$
      k_0\ge q_1=\rounddown{\frac1{\sqrt{2d}-\rounddown{\sqrt{2d}}}} \ .
   $$
   Since
   $$
      \liminf\set{\sqrt{2d}-\rounddown{\sqrt{2d}}\Bigwith
      d\ge 1,\ \sqrt{2d} \mbox{ irrational}} = 0 \ ,
   $$
   we can then choose $d$ in such a way that $k_0>N$, and hence
   $\ell_0>k_0\sqrt{2d}>N\sqrt{2d}$. This gives (b).
   It follows from a calculation that then, after possibly repeating the
   argument with a larger $N$, the inequalities in (a) are satisfied as well.
\endproof

\begin{remark}\rm
   Even though the formula for $\eps(X,L)$ in Theorem \ref{theorem asf}
   involves the solutions of a diophantine equation,
   the values for $\eps(X,L)$
   can be effectively computed in terms of $d$,
   since the solutions of Pell's equation can be obtained
   via continued fractions.
   In order to illustrate the situation, we include here a table
   providing the explicit (rounded) values
   of $k_0$, $\ell_0$,
   $\eps(X,L)$ and $\sqrt{2d}$ for $1\le d\le 30$ (see Table 1).
   We know from Theorem \ref{theorem asf} that the curve $C_0$
   computing $\eps(X,L)$ is of $L$-degree $2dk_0$ or $4dk_0$ and has
   a point of multiplicity $l_0$ or $2l_0$ respectively.
   Notice in particular how close $\eps(X,L)$ is to the theoretical
   upper bound in the cases $d=23$ and $d=29$.
   The curve $C_0$ has multiplicity
   $24335$ or $48670$ for $d=23$ and $19603$ or $39206$ for $d=29$.
   It does not come as a surprise then that it is hard to find $C_0$
   geometrically.
\end{remark}

\begin{table}
\begin{minipage}{\textwidth}
   \centering
   \footnotesize
   \begin{minipage}{0.5\textwidth}
      \begin{tabular}{rrrll}
	    $d$ & $k_0$ & $\ell_0$ & $\eps(X,L)$ & $\sqrt{2d}$ \\
	    \hline
             1&         2&         3&    1.333333333&    1.414213562 \\
             2&        --&        --&    2.000000000&    2.000000000 \\
             3&         2&         5&    2.400000000&    2.449489743 \\
             4&         1&         3&    2.666666667&    2.828427125 \\
             5&         6&        19&    3.157894737&    3.162277660 \\
             6&         2&         7&    3.428571429&    3.464101615 \\
             7&         4&        15&    3.733333333&    3.741657387 \\
             8&        --&        --&    4.000000000&    4.000000000 \\
             9&         4&        17&    4.235294118&    4.242640687 \\
            10&         2&         9&    4.444444444&    4.472135955 \\
            11&        42&       197&    4.690355330&    4.690415760 \\
            12&         1&         5&    4.800000000&    4.898979486 \\
            13&        10&        51&    5.098039216&    5.099019514 \\
            14&        24&       127&    5.291338583&    5.291502622 \\
            15&         2&        11&    5.454545455&    5.477225575
      \end{tabular}
   \end{minipage}%
   \nolinebreak%
   \begin{minipage}{0.5\textwidth}
      \begin{tabular}{rrrll}
	    $d$ & $k_0$ & $\ell_0$ & $\eps(X,L)$ & $\sqrt{2d}$ \\
	    \hline
            16&         3&        17&    5.647058824&    5.656854249 \\
            17&         6&        35&    5.828571429&    5.830951895 \\
            18&        --&        --&    6.000000000&    6.000000000 \\
            19&         6&        37&    6.162162162&    6.164414003 \\
            20&         3&        19&    6.315789474&    6.324555320 \\
            21&         2&        13&    6.461538462&    6.480740698 \\
            22&        30&       199&    6.633165829&    6.633249581 \\
            23&      3588&     24335&    6.782329977&    6.782329983 \\
            24&         1&         7&    6.857142857&    6.928203230 \\
            25&        14&        99&    7.070707071&    7.071067812 \\
            26&        90&       649&    7.211093991&    7.211102551 \\
            27&        66&       485&    7.348453608&    7.348469228 \\
            28&         2&        15&    7.466666667&    7.483314774 \\
            29&      2574&     19603&    7.615773096&  7.615773106 \\
            30&         4&        31&    7.741935484&    7.745966692
      \end{tabular}
   \end{minipage}

   \bigskip
   \begin{minipage}{0.7\textwidth}
      Table 1: The Seshadri constants $\eps(X,L)$ of abelian surfaces
      $(X,L)$ of type $(1,d)$ for $1\le d\le 30$ when $\NS(X)\isom\Z$.
   \end{minipage}

\end{minipage}
\end{table}


\section{The nef cone of an abelian surface}

   We start with some remarks on cones.
   Let $V$ be a real vector space and let $\Lambda$ be a lattice in $V$.
   A {\em cone} in $V$ is a subset $C\subset V$ such that $\R^+\cdot C\subset C$.
   It is convex if and only if $C+C\subset C$.
   A cone is called {\em polyhedral} if there are finitely many elements
   $\liste v1r\in V$ such that
   $$
      C=\sum_{i=1}^r \R^+\cdot v_i \ ,
   $$
   and it is said to be {\em rational polyhedral} if the generators
   $\liste v1r$ can be chosen within $\Lambda$.
   The dual $C^*\subset V^*$ of a cone $C$ is the cone
   $$
      C^*=\set{w\in V^*\with \bilin wv\ge 0\mbox{ for all }v\in C} \ .
   $$
   One has $C=C^{**}$ if and only if $C$ is closed and convex.
   $C$ is (rational) polyhedral if and only if $C^*$ is.
   Furthermore, by Gordon's Lemma,
   $C$ is rational polyhedral if and only if the semi-group
   $C\cap\Lambda$ is finitely generated.
   (See e.g.\ \cite[Theorem 14.1 and \S\S19,20]{Roc70}
   for the elementary properties of cones mentioned in this paragraph.)

   Consider now
   a smooth projective variety $X$.
   Via the intersection product, the real vector space
   $$
      \N_1(X)\eqdef\set{\mbox{1-cycles on $X$ modulo numerical
      equivalence}}\tensor\R
   $$
   is dual to the N\'eron-Severi vector space
   $\NSR(X)=\NS(X)\tensor\R$.
   As usual denote by $\NE(X)$ the {\em cone of curves} on $X$,
   i.e.\ the convex cone in
   $\N_1(X)$ generated by the effective 1-cycles.
   The dual cone of $\NE(X)$ is
   the nef cone 
   $$
      \Nef(X)
      =\set{\lambda\in\NSR(X)\with\lambda\xi\ge 0\mbox{
      for all }\xi\in\NE(X)} \ ,
   $$
   and the dual of $\Nef(X)$ is in turn
   the {\em closed cone of curves}
   $\bNE(X)$, i.e.\ the closure of $\NE(X)$ in $N_1(X)$, so
   $$
      \NE(X)^{**}=\Nef(X)^*=\bNE(X) \ .
   $$
   By the Cone Theorem \cite{Mor82} one knows that
   $\bNE(X)$, and hence $\Nef(X)$, is rational polyhedral
   whenever $c_1(X)$ is ample.
   If $c_1(X)$ is not ample, however,
   the structure of
   $\bNE(X)$ can be quite hard to determine
   and it will in general
   depend in a subtle way
   on the geometry of $X$ (cf.\ \cite[\S4]{CKM}).
   A good example for this phenomenon is
   \cite{Kov94} where
   $\bNE(X)$ is studied for $K3$ surfaces.

   Let now $L$ be
   an $\R$-line bundle, i.e.\ an element of
   $\Pic(X)\tensor_\Z\R$.
   It is {\em ample} if the conditions of the Nakai-Moishezon criterion
   hold for $L$.  By \cite{CamPet90} this is equivalent to requiring
   that (the numerical equivalence class of)
   $L$ belong to the interior of the nef cone
   $\Nef(X)$.
   Since we can take the pull-back of an $\R$-line bundle by a
   morphism, there is no problem to extend the
   definition of Seshadri constants to $\R$-line bundles $L$:
   \be
      \eps(L,x)&=&\sup\set{\eps>0\with f^*L-\eps E} \\
               &=&\inf\set{\frac{L\cdot C}{\mult_x C}\Biggwith
               C\subset X\mbox{ irreducible curve}} \ ,
   \ee
   where $f:\Bl_x(X)\to X$ denotes the blow-up at $x$ and
   $E=f\inverse(x)$. Again one has $\eps(L,x)>0$ for all $x\in X$,
   if $L$ is ample.

   As we will see,
   knowledge on the structure of the nef cone can be useful in the
   computation of Seshadri constants.
   Suppose for instance that $\Nef(X)$ 
   is polyhedral
   (or, equivalently, that $\bNE(X)$ is polyhedral),
   i.e.\
   \begin{equation}\label{def polyhedral}
      \Nef(X)=\sum_{i=1}^r\R_0^+\cdot[N_i]
   \end{equation}
   with $\R$-line bundles $N_i$ on $X$, so that
   if $L\in\Pic(X)$ is any ample line bundle, we have
   $$
      L\equiv\sum_{i=1}^r a_i N_i
   $$
   with real numbers $a_i > 0$.
   The Seshadri constant of $L$ at a point $x\in X$ is then clearly
   bounded below in terms of the numbers $a_i$ and the Seshadri
   constant of the line bundle $\sum_{i=1}^r N_i$:
   \begin{equation}\label{polyhedral bound}
      \eps(L,x)\ge\min\{a_1,\dots,a_r\}\cdot\eps\(\sum_{i=1}^r N_i, x\)
      \ .
   \end{equation}
   Note that $\sum_{i=1}^r N_i$ is ample, so that the bound is indeed
   non-trivial (and in many cases even sharp, as we will see below).

\begin{example}\rm
   Consider a principally polarized abelian surface $(X,L_0)$ with
   endomorphism ring $\End(X)\isom\Z[\sqrt{d}]$, where $d$ is a
   square-free positive integer.  Application of Shimura's theory
   shows that there is a two-dimensional family of such surfaces for
   any such $d$ (see \cite{Bir94}).
   We know that $\eps(L_0)=\frac 43$ (see \cite{Ste} or
   Theorem \ref{theorem asf}).  What can we say about the
   Seshadri constants of the other ample line bundles on $X$?
   First recall that the principal polarization induces
   an isomorphism of groups
   $$
      \NS(X)\to\End^s(X),\ L\mapsto\phi_{L_0}\inverse\circ\phi_L\ ,
   $$
   where $\End^s(X)\subset\End(X)$ is the subgroup of endomorphisms
   that are symmetric with respect to the Rosati involution
   $f\mapsto\phi_{L_0}\inverse\circ f\circ\phi_{L_0}$ on $\End(X)$.
   The endomorphism $\sqrt{d}$ has the characteristic polynomial
   $t^2-d$, hence the corresponding line bundle $\Lsd\in\NS(X)$
   satisfies $\Lsd^2=-2d$ and the classes of $L_0$ and $\Lsd$ yield an
   orthogonal (with respect to the intersection form) basis of
   $\NS(X)$.  By the 
   version of the Nakai-Moishezon Criterion given in
   \cite[Corollary 4.3.3]{LB}, a line bundle
   $$
      aL_0+b\Lsd \quad\mbox{ with } a,b\in\Z
   $$
   is ample if and only if
   $$
      (aL_0+b\Lsd)^2>0\quad\mbox{ and }\quad(aL_0+b\Lsd)L_0>0 \ .
   $$
   This implies that the nef cone is given by
   \be
      \Nef(X)&=&\set{xL_0+y\Lsd\in\R^2\Bigwith x\ge\sqrt{d}|y|} \\
             &=&\R_0^+(\sqrt{d}L_0+\Lsd)+\R_0^+(\sqrt{d}L_0-\Lsd) \ .
   \ee
   Note that it is polyhedral, but not rational polyhedral.
   For the generators $N^\pm=\sqrt{d}L_0\pm\Lsd$ we have
   $$
      \eps(N^++N^-)=\eps(2\sqrt dL_0)=2\sqrt d\cdot\eps(L_0)=\frac 83\sqrt
      d \ .
   $$
   so that \eqnref{polyhedral bound} gives us the lower bound
   $$
      \eps(aN^++bN^-)\ge\frac 83\sqrt d\cdot\min\{a,b\} \ .
   $$
   This bound is indeed sharp, as one sees by
   taking $a=b=\frac{1}{2\sqrt d}$, where one gets
   $aN^++bN^-=L_0$.
\end{example}

   The previous example shows that knowledge about the nef
   cone of a variety can be useful
   for the computation of Seshadri constants.
   Consider now the case when $X$ is an abelian variety.
   The most pleasant case, of course, is that $\Nef(X)$ is rational
   polyhedral, i.e.\ the case where the generators $N_i$ in
   \eqnref{def polyhedral} can be taken as (integral) line bundles
   in $\Pic(X)$.  By \cite{Bau1}, this happens if and
   only if $X$ is
   isogenous to a product
   $$
      X_1\times\dots\times X_r
   $$
   of mutually non-isogenous abelian varieties $X_i$ of Picard number
   one.  In general, the structure of $\Nef(X)$ can be more
   complicated.
   For the surface case,
   the following theorem gives the complete classification.
   We think that such a list is interesting, quite apart
   from its potential application on Seshadri constants.
   In the formulation of the proposition we distinguish the cases
   according to the rank of $\NS(X)$.
   Recall that $1\le\rank\NS(X)\le 4$ for every abelian surface $X$.

\begin{theorem}
   Let $X$ be an abelian surface, and let $\rho(X)=\rank\NS(X)$
   denote its Picard number.
   \begin{pitems}
   \item[(a)]
      If $\rho(X)=1$, then $\Nef(X)\isom\R_0^+$.
   \item[(b)]
      Suppose $\rho(X)=2$ and let $L,L'$ be line bundles whose classes
      generate $\NS(X)\tensor\Q$, with $L$ being ample.
      Consider the integer
      $$
         \delta(L,L')\eqdef (L\cdot L')^2-(L^2)((L')^2) \ .
      $$
      Then $\Nef(X)$ is polyhedral,
      $$
         \Nef(X)\isom\set{(x,y)\in\R^2\Bigwith
         x\ge\sqrt{\delta(L,L')}\cdot |y|}
         \ .
      $$

      If $\delta(L,L')$ is a square, then $\Nef(X)$ is rational
      polyhedral. In this case $X$ is isogenous to a product $E_1\times E_2$
      of non-isogenous elliptic curves $E_i$ with $\End(E_i)=\Z$.

      If $\delta(L,L')$ is not a square, then $\Nef(X)$ is irrational
      polyhedral. In this case $X$ is simple and has real or complex
      multiplication.

   \item[(c)]
      Suppose $\rho(X)=3$. Then $\Nef(X)$ is a cone over a circle :
      $$
         \Nef(X)\isom\set{(x,y,z)\in\R^3\with z^2\le xy,\ x\ge 0,\ y\ge 0} \ .
      $$
      Either $X$ is isogenous to the self-product $E\times E$ of an
      elliptic curve $E$ with $\End(E)\isom\Z$,
      or $X$ is simple and has indefinite quaternion multiplication
      (i.e.\ $\End_\Q(X)$ is an indefinite quaternion algebra).

   \item[(d)]
      Suppose $\rho(X)=4$. Then $X$ is isogenous to the self-product
      $E\times E$ of an elliptic curve $E$ with complex multiplication,
      and $\Nef(X)$ is a cone over a half-sphere:
      $$
         \Nef(X)\isom\set{(x,y,z,t)\in\R^4\with z^2+t^2\le xy,\ x\ge 0,\ y\ge 0} \ .
      $$
   \end{pitems}
\end{theorem}

\proof
   (a) is clear. As for (b), consider the line bundle
   $$
      M \eqdef (L\cdot L')L-L^2\cdot L'\in\Pic(X) \ .
   $$
   We have $L\cdot M=0$ and $M^2=-\delta(L,L')\cdot L^2$.
   Note that $\delta(L,L')>0$, i.e.\ $M^2<0$, by the Hodge index theorem.
   By assumption, every line bundle
   on $X$ is algebraically equivalent to a bundle $aL+bM$ with suitable
   $a,b\in\Q$. Now, $aL+bM$ is ample if and only if
   $$
      (aL+bM)^2>0 \and (aL+bM)L>0 \ ,
   $$
   and these conditions are satisfied if and only if
   $$
      a^2 > b^2\cdot\delta(L,L') \ .
   $$
   This implies the statement about the nef cone and shows that it is
   in any event polyhedral.
   The quadratic form
   $$
      \psi:\Q^2\to\Q,\quad (a,b)\mapsto (aL+bM)^2
   $$
   represents zero (non-trivially) if and only if $\sqrt{\delta(L,L')}$
   is rational,
   i.e.\ if and only if $\Nef(X)$ is {\em rational} polyhedral.
   Now, if $\psi$ represents zero,
   $$
      0=\psi(a,b)=\(a^2-b^2\cdot\delta(L,L')\)L^2 \ ,
   $$
   then $(aL+bM)L>0$,
   if we choose $a>0$.
   But this implies that the class $aL+bM$ is effective
   and is represented by a multiple of an elliptic curve. So $X$ is isogenous
   to a product $E_1\times E_2$ of elliptic curves in this case.
   The condition $\rho(X)=2$ implies that we must have $\End(E_i)=\Z$ and
   that the $E_i$ are non-isogenous.
   In the alternative case, i.e.\ when $\psi$ does not represent zero,
   then $X$ is simple, and the classification of endomorphisms algebras
   of simple abelian varieties (see \cite[Chap.\ 5]{LB} or
   \cite[Sect.\ 21]{Mum74a})
   shows that then $X$ has either real or
   complex multiplication.

   (c) Consider first the case when $X$ is simple.
   One sees from the classification of endomorphism algebras that
   then the assumption $\rho(X)=3$ implies
   that $X$ has indefinite quaternion multiplication. There is an
   isomorphism of the algebra $\End_\R(X)$ with $M_2(\R)$
   under which the Rosati involution
   on $\End_\R(X)$ corresponds to matrix transposition. The ample classes
   in $\NS_\R(X)$ correspond under the composed map
   $$
      \NS_\R(X)\isomto\End^s_\R(X)\isomto\Sym_2(\R)
   $$
   to the positive definite matrices.
   This implies the statement on the nef cone.
   Suppose now that $X$ is non-simple.  The condition $\rho(X)=3$
   implies that then $X$ is isogenous to a product $E\times E$, where
   $E$ is an elliptic curve with $\End(E)=\Z$.
   The nef cone of $E\times E$, and hence also the nef cone of $X$,
   can be described explicitly in terms of a suitable basis
   of $\NS_\Q(E\times E)$.  Taking $\Delta\subset E\times E$ to be the
   diagonal, the classes of $E_1=E\times 0, E_2=0\times E$ and
   $F=\Delta-E_1-E_2$ generate $\NS_\Q(E\times E)$.
   A line bundle
   $$
      M(a_1,a_2,b)\eqdef a_1E_1+a_2E_2+bF
   $$
   is ample if and only if
   $$
      M(a_1,a_2,b)^2>0 \and M(a_1,a_2,b)\cdot E_i>0 \quad\mbox{ for } i=1,2
      \ .
   $$
   But these conditions are equivalent to $b^2<a_1a_2$
   and $a_i>0$ respectively, and this proves
   the assertion.

   (d) The condition $\rho(X)=4$ implies that $X$ is isogenous to $E\times E$,
   where $E$ is an elliptic curve with complex multiplication
   (see \cite{ShiMit74}).

   Then
   $$
      \End_\Q(X)=M_2(\End_\Q(E))=M_2(\Q(\sqrt d)) \ ,
   $$
   with an integer $d<0$.
   Taking the Rosati involution with respect to a product polarization,
   we have
   $$
      \End^s(X)
      =\set{\matr{f_1}{f_2}{f_3}{f_4} \Biggwith f_1,f_4\in\End^s(E)
      \mbox{ and } f_2'=f_3 }
   $$
   The ample classes in $\NS_\R(X)$ then correspond to the positive definite
   matrices of the form
   $$
      \matr{\alpha}{\beta+\gamma\sqrt d}{\beta-\gamma\sqrt d}{\eta}
   $$
   with $\alpha,\beta,\gamma,\eta\in\R$, and from this follows the assertion
   on the nef cone of $X$.
\endproof


\section{Multiple point Seshadri constants on abelian surfaces}
\label{section multiple}

   Consider a smooth projective variety $X$ and an ample line bundle
   $L$ on $X$.  So far we have considered the Seshadri constant
   $\eps(L,x)$ of $L$ at a single point $x$ of $X$.  There is a
   natural generalization of this idea, a {\em multiple point Seshadri
   constant}, accounting for the positivity of $L$ along a finite set
   of points.
   It is quite obvious what the natural definition is:
   For distinct points $\liste x1k$ in $X$ one puts
   $$
      \eps(L,\liste x1k) \eqdef \sup\set{\eps\in\R\Biggwith
      f^*L-\eps\sum_{i=1}^k E_i\mbox{ nef}}
   $$
   where $f:Y\to X$ is the blow-up of $X$ at $\liste x1k$ and
   $E_i=f\inverse(x_i)$ for $1\le i\le k$ (cf.\ \cite[(5.16)]{Laz93}).
   As in the one-point case, one has the equivalent definition
   $$
      \eps(L,\liste x1k)=\inf_C \frac{L\cdot
      C}{\sum_{i=1}^k\mult_{x_i} C} \ ,
   $$
   where the infimum is taken over all irreducible curves $C\subset X$
   passing through at least one of the points $\liste x1k$.

   There are the upper and lower bounds
   \begin{equation}\label{elem multiple bounds}
      \frac 1k\min_{1\le i\le k}\eps(L,x_i)\le\eps(L,\liste
      x1k)\le\sqrt[n]{\frac{L^n}{k}} \ ,
   \end{equation}
   where $n=\dim X$.
   In fact, letting $\eps=\eps(L,\liste x1k)$, we have
   $$
      L^n-k\eps^n=\(f^*L-\eps\sum_{i=1}^k E_i\)^{\!\!n}\ge 0 \ ,
   $$
   since $f^*L-\eps\sum E_i$ is nef, and this implies the second
   inequality in \eqnref{elem multiple bounds}.  Further, putting
   $\delta=\min\set{\eps(L,x_i)\with 1\le i\le k}$, the line bundle
   $$
      k\cdot f^*L-\delta\sum_{i=1}^k E_i=\sum_{i=1}^k(f^*L-\delta E_i)
   $$
   is nef, and hence $\eps(L,\liste x1k)\ge\delta/k$, which gives the
   first inequality in \eqnref{elem multiple bounds}.

   While already one-point Seshadri constants are subtle invariants,
   multiple point Seshadri constants are extremely hard to control.
   Suffice it to say that Nagata's famous conjecture relating the degrees
   and multiplicities of curves in the projective plane
   can equivalently
   be formulated as a statement on multiple point
   Seshadri constants on $\P^2$:

\begin{varthm}{Nagata's Conjecture}[cf.\ \cite{Nag59}]
   For general points $\liste x1k\in\P^2$ with $k\ge 9$ one has
   $$
      \eps(\O_{\P^2}(1), \liste x1k)=\frac 1{\sqrt k} \ .
   $$
\end{varthm}

   In other words, the conjecture says that the multiple point
   Seshadri constant of $\O_{\P^2}(1)$ at $\ge 9$
   general points should have its
   maximal possible value.  This is known to be true whenever $k$ is a
   square.

   In light of these facts
   it seems hardly surprising that only few general results on multiple
   point Seshadri constants are available (e.g.\  \cite{Kue96},
   \cite{Xu94}).  Let us consider here the case of abelian surfaces.
   By homogeneity,
   the number $\eps(L,\liste x1k)$ depends then only on the
   differences $x_i-x_1$, $2\le i\le k$, and we have
   $$
      \eps(L,\liste x1k)\ge \frac 1k\eps(L,x_1)=\frac 1k\eps(L) \ .
   $$
   The following result shows that in general (i.e.\
   if either $(X,L)$ is general or if
   the points $\liste x1k$ are general on $X$) one has the
   strict inequality whenever $k\ge 2$.
   In fact, we show that equality can only hold
   for trivial reasons:

\begin{proposition}\label{single multiple sesh}
   Let $X$ be an abelian surface and let $L$ be an ample line bundle
   on $X$.
   Suppose that $\liste x1k\in X$
   are distinct points with $k\ge 2$.
   If
   $$
      \eps(L,\liste x1k)=\frac 1k \eps(L) \ ,
   $$
   then $X$ contains an elliptic curve $E$ with the property
   $$
      L\cdot E = k\cdot\eps(L,\liste x1k)
   $$
   and all the points $\liste x1k$ lie on $E$.
\end{proposition}

   The proposition will be deduced from the following result, which
   provides a lower bound on the multiple point Seshadri constant in
   terms of $L^2$ and $k$.

\begin{proposition}\label{bound multiple sesh}
   Let $X$ be an abelian surface and let $L$ be an ample line bundle.
   Then for every choice of points $\liste x1k\in X$ with $k\ge 1$ we are
   in one of the following two cases:
   \begin{items}
   \item[(a)]
      $$
         \eps(L,\liste x1k)\ge
         \begin{casearray}
            \displaystyle\frac{\sqrt{2L^2}}k & \mbox{ , if } k\ge 4
            \\[1em]
            \displaystyle\frac{\sqrt{L^2}}{2\sqrt 2}\sqrt{\frac{8-k}k} &
            \mbox{ , if } 1\le k\le 3
         \end{casearray}
      $$
      or
   \item[(b)]
      $X$ contains an elliptic curve $E$ such that
      $$
         \eps(L,\liste x1k)=\frac{L\cdot E}{\#\set{i\with x_i\in E}} \ .
      $$
   \end{items}
\end{proposition}

   As one might expect, the bound in (a) increases
   with $L^2$ and decreases with $k$.  It differs from the theoretical
   upper bound by a factor of order $1/\sqrt k$.
   Note that the appearance of a case dealing with the potential
   influence of elliptic curves on the Seshadri constant is quite
   inevitable: Polarized abelian surfaces $(A,L)$
   can contain elliptic curves of any
   given $L$-degree, no matter how large $L^2$ is.

\proofof{Proposition \ref{bound multiple sesh}}
   We consider first the contribution of the non-elliptic curves
   to the Seshadri constant, i.e.\ the number
   $$
      \eps'(L,\liste x1k)\eqdef\inf\set{\frac{L\cdot C}{\sum_{i=1}^k
      \mult_{x_i} C}\Biggwith C\subset X
      \mbox{ ample irreducible curve}} \ .
   $$
   Note that the self-intersection of the curves $C\subset X$ in
   question
   is bounded below by
   \begin{equation}\label{multiple mult estimate}
      C^2\ge 2+\sum_{i=1}^k\mult_{x_i}C(\mult_{x_i}C-1) \ .
   \end{equation}
   In fact, since $X$ contains no rational curves and since all
   elliptic curves are smooth, we have $p_g(C)\ge 2$, and hence
   $$
      \frac 12 C^2-1\ge p_a(C)-p_g(C)\ge\sum_{i=1}^k
      {\mult_{x_i}C\choose 2 } \ ,
   $$
   which implies \eqnref{multiple mult estimate}.

   Now, since $\sum_{i=1}^k(\mult_{x_i}C)^2\ge (1/k)(\sum_{i_1}^k
   \mult_{x_i}C)^2$, the inequality \eqnref{multiple mult estimate} gives
   a quadratic relation for the sum $\sum_{i=1}^k\mult_{x_i} C$,
   $$
      \(\sum_{i=1}^k\mult_{x_i}C\)^2 - k\sum_{i=1}^k\mult_{x_i}C +
      k(2-C^2) \le 0 \ ,
   $$
   which tells us that
   $$
      \sum_{i=1}^k\mult_{x_i}C\le \frac k2+\sqrt{k\(\frac k4+C^2-2\)} \ .
   $$
   Using now the Hodge index theorem for the line bundles $L$ and
   $\O_X(C)$, this bound on the multiplicities yields a bound on
   the number $\eps'(L,\liste x1k)$,
   $$
      \eps'(L,\liste x1k)\ge\inf_C
      \frac{L\cdot C}{\frac k2+\sqrt{k\(\frac k4+\frac{(L\cdot
      C)^2}{L^2}-2\)}} \ ,
   $$
   where the infimum is taken over the ample irreducible curves $C\subset X$.
   Consider now for fixed numbers $k$ and $L^2$ the real-valued
   function
   $$
      f:t\mapsto\frac t{\frac k2+\sqrt{k(\frac k4+\frac{t^2}{L^2}-2)}}
      \ .
   $$
   For $k<8$ its minimum lies at the point
   $$
      t_0=\sqrt{2L^2}\sqrt{\frac{8-k}k} \ ,
   $$
   whereas for $k\ge 8$ the function is increasing.  Note next that
   $$
      t=L\cdot C\ge\sqrt{L^2}\sqrt{C^2}\ge \sqrt{2L^2}\eqdef t_1
   $$
   and that $t_0\le t_1$ for $k\ge 4$.  We conclude that
   $$ 
      \eps'(L,\liste x1k)\ge\min f|[t_1,\infty) = f(t_1) =   
      \frac{\sqrt{2L^2}}k
   $$
   for $k\ge 4$, and
   $$
      \eps'(L,\liste x1k)\ge f(t_0)=\frac{\sqrt{L^2}}{2\sqrt 2}
      \sqrt{\frac{8-k}k}
   $$
   for $k\le 3$.

   Finally, if $\eps(L,\liste x1k)<\eps'(L,\liste x1k)$, then the
   multiple point Seshadri constant must be computed by an elliptic
   curve $E$, and one has
   $\eps(L,\liste x1k)=L\cdot E/(\#\set{i\with x_i\in E})$.
\endproof

   We now give the

\proofof{Proposition \ref{single multiple sesh}}
   We will first show that under the hypothesis of the proposition
   the Seshadri constant $\eps(L,\liste x1k)$
   must be computed by an elliptic curve.  Suppose to the contrary that
   $\eps(L,\liste x1k)=\eps'(L,\liste x1k)$ (in the notation of the
   proof of the previous proposition).
   Then, using Proposition
   \ref{bound multiple sesh} and the hypothesis we get
   $$
      \eps(L,x_1)=k\cdot\eps(L,\liste x1k)\ge
         \begin{casearray}
            {\sqrt{2L^2}} & \mbox{ , if } k\ge 4 \\
            \frac{1}{2\sqrt 2}\sqrt{k(8-k)L^2} & \mbox{ , if }
            1\le k\le 3
         \end{casearray}
   $$
   On the other hand, one always has $\eps(L,x_1)\le\sqrt{L^2}$,
   so that
   $$
      \frac 1k\sqrt{L^2}\ge
         \begin{casearray}
            \frac{\sqrt{2L^2}}k & \mbox{ , if } k\ge 4 \\
            \frac{\sqrt{L^2}}{2\sqrt 2}\sqrt{\frac{8-k}k} & \mbox{ , if }
            1\le k\le 3
         \end{casearray}
   $$
   In case $k\ge 4$ we certainly have a contradiction, and for $1\le
   k\le 3$ we get
   $(1/k)\sqrt{2d}\ge(1/2)\sqrt{2d(8-k)/(2k)}$,
   which is impossible due to the assumption $k\ge 2$.

   So we conclude that $\eps(L,\liste x1k)$ is computed by an elliptic
   curve $E$, and it remains to show that all the points $\liste x1k$ lie
   on $E$.  But, putting $\ell=\#\set{i\with x_i\in E}$, we find
   $$
      \frac{L\cdot E}{\ell}=\eps(L,\liste x1k)=\frac 1k\eps(L,x_1)
      \le \frac{L\cdot E}k \ ,
   $$
   which implies $\ell=k$, and this completes the proof.
\endproof



\bigskip
   Mathematisches Institut,
   Uni\-ver\-si\-t\"at Er\-lan\-gen-N\"urn\-berg,
   Bis\-marck\-stra{\ss}e $1\frac12$,
   D-91054 Erlangen,
   Germany.
   {\em E-mail address:} {\tt bauerth@mi.uni-erlangen.de}


\end{document}